\documentclass[a4paper]{article}
\usepackage{amsmath,amsthm,amsfonts,thmtools}
\usepackage{amssymb}
\usepackage[cmintegrals]{newtxmath}
\usepackage{graphicx}
\usepackage{listings}
\usepackage[font=footnotesize]{caption}
\usepackage[font=scriptsize]{subcaption}
\usepackage{euscript}
\usepackage{units}
\usepackage{enumerate}
\usepackage{enumitem}
\usepackage{breqn}
\usepackage{xcolor}
\usepackage{todonotes}
\usepackage{booktabs}
\usepackage{tikz}
\usepackage{comment}
 
\DeclareMathAlphabet{\pazocal}{OMS}{zplm}{m}{n}

\usepackage[bookmarks=true,hidelinks]{hyperref}
\usepackage{bbm}
\usepackage{multirow}
\usepackage{algpseudocode}  % for the algorithmic environment
\usepackage{mathtools}
\usepackage{subcaption} 
\usepackage{graphicx}
\usepackage{pgfplots}
\usepackage{float}
\usepackage{algpseudocode}
\usepackage{bbold}
\usepackage{cite}
\usepackage{microtype}
\usepackage{rotating}
\usepackage{enumitem}
\usepackage{bm}
\usepackage{hyperref}
\providecommand{\keywords}[1]
{
  \small	
  \textbf{\textit{Keywords:}} #1
}
\numberwithin{equation}{section}
\newtheorem{theorem}{Theorem}[section]

\newtheorem{algorithm}[theorem]{Algorithm}

\newtheorem{appendixthm}{}[section] 
\numberwithin{equation}{section}

\theoremstyle{appendixstyle}

\newtheoremstyle{appendixstyle} % Custom style name
  {3pt}   % Space above
  {3pt}   % Space below
  {\normalfont} % Body font
  {}      % Indent amount
  {\bfseries} % Theorem head font
  {.}     % Punctuation after theorem head
  {.5em}  % Space after theorem head
  {\thmnumber{#2}} % Custom head specification, just show the number

\theoremstyle{definition}

\newtheoremstyle{myremarkstyle}{}{}{}{}{\bfseries}{.}{ }{}
\theoremstyle{myremarkstyle}
\declaretheorem[name=Remark,qed=$\blacksquare$,numberlike=theorem]{remark}

\makeatletter
\newcommand*{\intavg}{%
  \mint@l{-}{}%
}
\newcommand*{\mint@l}[2]{%
  \@ifnextchar\limits{%
    \mint@l{#1}%
  }{%
    \@ifnextchar\nolimits{%
      \mint@l{#1}%
    }{%
      \@ifnextchar\displaylimits{%
        \mint@l{#1}%
      }{%
        \mint@s{#2}{#1}%
      }%
    }%
  }%
}
\newcommand*{\mint@s}[2]{%
  \@ifnextchar_{%
    \mint@sub{#1}{#2}%
  }{%
    \@ifnextchar^{%
      \mint@sup{#1}{#2}%
    }{%
      \mint@{#1}{#2}{}{}%
    }%
  }%
}
\def\mint@sub#1#2_#3{%
  \@ifnextchar^{%
    \mint@sub@sup{#1}{#2}{#3}%
  }{%
    \mint@{#1}{#2}{#3}{}%
  }%
}
\def\mint@sup#1#2^#3{%
  \@ifnextchar_{%
    \mint@sub@sup{#1}{#2}{#3}%
  }{%
    \mint@{#1}{#2}{}{#3}%
  }%
}
\def\mint@sub@sup#1#2#3^#4{%
  \mint@{#1}{#2}{#3}{#4}%
}
\def\mint@sup@sub#1#2#3_#4{%
  \mint@{#1}{#2}{#4}{#3}%
}
\newcommand*{\mint@}[4]{%
  \mathop{}%
  \mkern-\thinmuskip
  \mathchoice{%
    \mint@@{#1}{#2}{#3}{#4}%
        \displaystyle\textstyle\scriptstyle
  }{%
    \mint@@{#1}{#2}{#3}{#4}%
        \textstyle\scriptstyle\scriptstyle
  }{%
    \mint@@{#1}{#2}{#3}{#4}%
        \scriptstyle\scriptscriptstyle\scriptscriptstyle
  }{%
    \mint@@{#1}{#2}{#3}{#4}%
        \scriptscriptstyle\scriptscriptstyle\scriptscriptstyle
  }%
  \mkern-\thinmuskip
  \int#1%
  \ifx\\#3\\\else_{#3}\fi
  \ifx\\#4\\\else^{#4}\fi  
}
\newcommand*{\mint@@}[7]{%
  \begingroup
    \sbox0{$#5\int\m@th$}%
    \sbox2{$#5\int_{}\m@th$}%
    \dimen2=\wd0 %
    \let\mint@limits=#1\relax
    \ifx\mint@limits\relax
      \sbox4{$#5\int_{\kern1sp}^{\kern1sp}\m@th$}%
      \ifdim\wd4>\wd2 %
        \let\mint@limits=\nolimits
      \else
        \let\mint@limits=\limits
      \fi
    \fi
    \ifx\mint@limits\displaylimits
      \ifx#5\displaystyle
        \let\mint@limits=\limits
      \fi
    \fi
    \ifx\mint@limits\limits
      \sbox0{$#7#3\m@th$}%
      \sbox2{$#7#4\m@th$}%
      \ifdim\wd0>\dimen2 %
        \dimen2=\wd0 %
      \fi
      \ifdim\wd2>\dimen2 %
        \dimen2=\wd2 %
      \fi
    \fi
    \rlap{%
      $#5%
        \vcenter{%
          \hbox to\dimen2{%
            \hss
            $#6{#2}\m@th$%
            \hss
          }%
        }%
      $%
    }%
  \endgroup
}

\def\XXint#1#2#3{{\setbox0=\hbox{$#1{#2#3}{\int}$ }
		\vcenter{\hbox{$#2#3$ }}\kern-.6\wd0}}

\renewcommand{\geq}{\geqslant}
\renewcommand{\leq}{\leqslant}

\renewcommand{\epsilon}{\varepsilon}
\renewcommand{\phi}{\varphi}

\begin{document}

\title{Physics informed neural network for forward and inverse radiation heat transfer in graded-index medium.}

\author{K. Murari \footnotemark[1] 
	\and S. Sundar \footnotemark[1]
}

\date{\today}

\maketitle
\medskip
\centerline{$^*$Centre for Computational Mathematics and Data Science}
\centerline{Department of Mathematics, IIT Madras, Chennai 600036, India}
\centerline{kmurari2712@gmail.com, slnt@iitm.ac.in}

\begin{abstract}
Radiation heat transfer in a graded-index medium often suffers accuracy problems due to the gradual changes in the refractive index. The finite element method,  meshfree, and other numerical methods often struggle to maintain accuracy when applied to this medium. To address this issue, we apply physics-informed neural networks (PINNs)-based machine learning algorithms to simulate forward and inverse problems for this medium. We also provide the theoretical upper bounds. This theoretical framework is validated through numerical experiments of predefined and newly developed models that demonstrate the accuracy and robustness of the algorithms in solving radiation transport problems in the medium. The simulations show that the novel algorithm goes on with numerical stability and effectively mitigates oscillatory errors, even in cases with more pronounced variations in the refractive index. 
\end{abstract}

\medskip
%---------------------------------------------------------------


\section{Introduction}\label{sec:intro}
\noindent
Radiation heat transfer in graded index materials involves the transmission of thermal energy through electromagnetic waves within a material exhibiting varying refractive index. Typically, this index changes gradually from the center to the material's surface, impacting radiation propagation and heat transfer mechanisms. Understanding the complexities of heat transfer in the materials is paramount. Unlike homogeneous materials, which maintain a constant refractive index, graded index materials undergo gradual changes, significantly affecting thermal radiation propagation. The pivotal role of radiation across various applications justifies this focus, ranging from interpreting spectroscopic emissions in celestial bodies \cite{fokou2021radiation}, nuclear engineering \cite{larsen1986one}, biological tissues \cite{randrianalisoa2014effects}, thermal insulation \cite{baillis2013thermal}, gas turbines \cite{seaid2004efficient}, combustion systems \cite{ebrahimi2013zonal} and radiotherapy dose simulation \cite{pichard2017approximation}.
Numerous researchers have focused on analyzing radiation heat transmission in such media. Various numerical methods, including the Monte Carlo Method \cite{wei2019reverse}, Least-Squares Finite Element Method (LSFEM) \cite{liu2007least}, Discontinuous Finite Element Method (DFEM) \cite{feng2018discontinuous}, Finite Volume Method (FVM) \cite{liu2006finite}, Petrov-Galerkin Method (MLPG) \cite{liu2006meshless}, Least-Squares Spectral Element Method (LSSEM) \cite{zhao2012deficiency}, and Spectral Element Method (SEM) \cite{wang2012meshless}, have been employed to simulate radiation transport in these complex systems. However, traditional methods like LSFEM, the Galerkin Finite Element Method (GFEM) \cite{zhao2012deficiency}, and specific mesh-free methods \cite{zhao2013second} suffer from various deficiencies. These methods often exhibit noticeable high-frequency errors or unwanted fluctuations, despite generally providing results that are close to the true solution. Additionally, the Lattice Boltzmann Method (LBM) \cite{liu2020solving}, \cite{yi2016lattice}, and Multi-Relaxation Time (MRT) Lattice Boltzmann Method \cite{feng2021performance}, have been applied, but challenges persist in achieving robust solutions.\\

In recent years, deep learning is an alternative way to avoid the curse of dimensionality. Due to this reason, deep learning has recently been an essential tool in modern technology and advanced research in the last few years. This learning contains many layers of transformations and scalar nonlinearities. Deep learning techniques demonstrate strong capabilities in approximating highly nonlinear functions.  Their computational framework, which employs statistical learning and large-scale optimization methods in conjunction with contemporary hardware and software, enhances their capacity to address nonlinear and high-dimensional partial differential equations (PDEs).
\subsection*{Nomenclature}
\begin{center}
\setlength{\fboxrule}{0.5pt}
\setlength{\fboxsep}{5pt}
\fbox{
\begin{minipage}{0.95\textwidth}
\begin{tabular}{p{0.1\textwidth} p{0.35\textwidth} p{0.1\textwidth} p{0.35\textwidth}}  % Adjust column widths
    $c_{o}$     & speed of light                & $ \sigma $    & activation function \\
    $n$         & refractive index              & $T$         & temperature distribution [K] \\
    $s$         & position vector               & $\mathbf{1}_{[.]}$ & unit step function\\
    $k_e$       & extinction coefficient        & $\rho$      & diffuse reflectivity \\
    $\boldsymbol\Omega$    & direction vector              & $\epsilon_{\bar{\omega}}$ & emissivity \\
    $\Phi$      & scattering phase function     & $\Theta$    & dimensionless temperature \\
    $k_a$       & absorption coefficient        &  $ L$ &      length \\
    $k_s$       & scattering coefficient        & $\mu$       & directional cosine \\
    $T_{g}$     & medium temperature [K]        & $\rho_{m}$  & saterial density  \\
    $\Omega'$   & incoming direction            & $\sigma_B$  & Stefan-Boltzmann constant \\
    $\hat{\boldsymbol{n}}_{\boldsymbol{\omega}}$ & normal vector to surface & $I_{c}$     & collimated intensity from laser \\
    $\Omega_{rs}$ & incident direction of specular reflection & $I_{d}$ & diffuse intensity \\
    $I_{bw}$    & black body intensity          & $ \boldsymbol{\omega}$ & scattering albedo \\
    $I$         & radiative intensity           & $b$         & black-body state \\
    $S$         & source term of differential equation & $0, L$ & left and right boundaries \\
    $G(s,t)$    & incident radiation            & $\mathbf{i},\mathbf{j},\mathbf{k}$ & unit vectors \\
    $q(s,t)$    & heat flux [W m$^{-2}$]        & $\lambda$   & learning rate \\
    $\rho^{spce}$ & specular reflectivity        & $\lambda_{reg}$ & regularization parameters \\
    $K-1$       & number of hidden layers       &   $\Omega$    & outgoing direction          \\
    $n_{\theta}$ & number of times model retrained in parallel  & $tanh$ & hyperbolic tangent(activation function)\\
    $ N_{int} $ & number of interior points & $N_{sb}$ & number of Sobol points\\
    $N_{tb}$ & number of temporal points & $N_{d}$  & number of data points \\
    $N$ & number of collocations or training points &  & \\
\end{tabular}
\end{minipage}
}
\end{center}
Deep neural networks (DNNs) are well known for their ability to universally approximate functions, a property that holds under certain conditions, as demonstrated by Cybenko \cite{cybenko1989approximations}, Hornik et al. \cite{hornik1989multilayer}, Barron \cite{barron1993universal}, and Yarotsky \cite{yarotsky2017error}. This property makes them suitable as trial functions for solving PDEs, often by minimizing the residual of the PDE at selected collocation points within the domain. The application of deep learning using DNNs has revolutionized numerous domains, including image and text categorization, machine vision, computational linguistics, voice recognition, self-governing systems, robotics, artificial intelligence in gaming, medical analysis, pharmaceutical research, climate simulation, financial prediction, and protein structure determination. These advancements highlight the versatility and effectiveness of DNNs in tackling complex problems across diverse domains.
 % This algorithm has been applied in supervised, semi-supervised, and unsupervised forms. It provides solutions for forward and inverse modeling within a unified optimization framework. PINN showcases its adaptability by accommodating both continuous and discrete forms of PDEs for both forward and inverse modeling endeavors. Continuous time-dependent model approximation relies on spatio-temporal function estimation. On the other hand, when dealing with discrete time modeling, it tackles exact implicit Runge-Kutta time-stepping techniques, supporting an unrestricted number of time increments. This novel algorithm has extended into XPINN\cite{jagtap2020extended}, cPINN\cite{jagtap2020conservative}, Parallel PINN\cite{shukla2021parallel}, and Gradient-enhenced\cite{yu2022gradient}  and authors also provided a Deepxde\cite{lu2021deepxde} library to solve PDE with PINN. \cite{dolean2022finite} and \cite{moseley2023finite} modified this algorithm with domain decomposition approch. \cite{de2022error} and \cite{de2024error} gives the theoretical error bound of XPINN\cite{jagtap2020extended} for navier stocks equation and kolmogorov pdes respectively. More importantly, there still needs to be more understanding of why these models occasionally fail to train. \cite{wang2022and} explore these issues by examining them through the framework of the neural tangent kernel. \\
A technique involves utilizing DNNs founded on explicit or partially explicit representation formulas applicable to parabolic and elliptic partial differential equations (PDEs). Researchers leverage this compositional structure to enhance the approximation capabilities of DNNs. Researchers such as \cite{hornik1989multilayer}, \cite{evans2018novo}, and \cite{beck2021solving}, among others, have introduced and examined this technique for a range of parametric elliptic and parabolic PDEs. In \cite{lagaris1998artificial}, authors discussed a similar approach for approximating linear transport equations with DNNs.
A key deep learning algorithm, Physics-Informed Neural Networks (PINNs) \cite{raissi2019physics} is a mesh-free method. PINNs has been successfully applied in supervised, semi-supervised, and unsupervised learning frameworks. PINNs provide solutions for both forward and inverse modeling problems within a unified optimization framework. This algorithm is remarkably adaptable to both continuous and discrete forms of PDEs. For continuous time modeling, PINNs approximate solutions by estimating spatio-temporal functions. For discrete time modeling, PINNs utilize implicit Runge-Kutta time-stepping techniques, allowing for an arbitrary number of time steps. This flexibility has led to the development of extensions such as XPINN \cite{jagtap2020extended}, cPINN \cite{jagtap2020conservative}, Parallel PINN \cite{shukla2021parallel}, and Gradient-Enhanced PINN \cite{yu2022gradient}. In addition, the Deepxde library \cite{lu2021deepxde} facilitates solving PDEs using PINNs. Recent studies, such as those by Dolean et al. \cite{dolean2022finite} and Moseley et al. \cite{moseley2023finite}, have applied domain decomposition approaches to modify PINN, and theoretical error bounds for XPINN have been established for the Navier-Stokes equations \cite{de2022error} and Kolmogorov PDEs \cite{de2024error}. However, there is still limited understanding of why these models occasionally fail to train effectively. Wang et al. \cite{wang2022and} explore this issue through the framework of the neural tangent kernel. A recent study by Mishra and Molinaro examines the generalization error of addressing forward problems \cite{mishra2023estimates} and inverse problems \cite{mishra2022estimates} across various linear and nonlinear partial differential equations (PDEs). The authors have also estimated generalized error bounds for the problems. The authors also worked on nonlinear dispersive PDEs \cite{bai2021physics}, where they derived generalized error bounds for forward problems. The study investigates the stability characteristics of the underlying PDE, utilizing these features to assess generalization errors in connection with training errors. Furthermore, \cite{de2024wpinns} introduces an innovative modification of PINNs called weak PINNs (wPINNs).\\

PINNs have been used to simulate both forward and inverse problems in the context of radiation transport. Mishra et al. \cite{mishra2021physics} developed a PINN-based algorithm for simulating radiation transport, establishing generalized error bounds for forward problems. Other studies, such as those by Huhn et al. \cite{huhn2023physics} and Riganti et al. \cite{riganti2023auxiliary}, have further applied PINNs to solve radiation transport problems in various configurations. The present work simulates both new and existing models for radiation heat transfer in graded-index media. We employ unsupervised PINN algorithms to simulate both predefined and newly developed models, demonstrating improved performance over traditional methods such as finite element methods (e.g., GFEM, LSFEM, DFEM), mesh-free methods, and MRT lattice Boltzmann methods. The results show that our PINN-based approach mitigates the challenges of oscillatory errors and achieves stable, accurate results even with complex variations in the refractive index.

The paper is structured as follows: Section \ref{sec:2} outlines the mathematical model and methodology, including the problem statement, PINN approximation, and key components such as the model, domain, quadrature rules, neural networks, residuals, loss functions, and optimization. Section  \ref{sec:3} presents numerical simulations and accuracy verification of the proposed method, while Section  \ref{sec:4} discusses the conclusion. An appendix is provided to estimate generalization errors and analyze the impact of optical parameters on radiometric quantities.
\section{Problem statement and PINN approximation}\label{sec:2}
Precise prediction of radiative transfer in these circumstances relies on resolving the Radiative Transfer Equation (RTE), an intricate integro -differential equation encompassing seven variables: three spatial coordinates, polar and azimuthal angles, time, and spectral dimensions. The complexity increases in graded-index media, where the curved radiation path, described by the Fermat principle, introduces an additional layer of intricacy, rendering analytical solutions unattainable except in specific limiting cases. Therefore, developing accurate, simple, and efficient tools to solve the RTE in these scenarios becomes imperative for advancing understanding and practical applications. Traditional methods, such as mesh-free or finite-element approaches, struggle with these complexities due to the difficulty of accurately capturing the light's curved path in graded-index media, which often leads to numerical instability and high computational cost. Snell's Law and the Fermat principle \cite{liu2006meshless} are closely related concepts that describe different aspects of light propagation but are fundamentally connected. Snell's Law \cite{arbab2016compton} quantitatively describes how light bends or refracts when transitioning between different optical mediums with varying refractive indices. It defines a mathematical connection between the angles of incidence and refraction.
In the context of Snell's Law, the light path between two points is the one that reduces the time needed for light to travel from the source to the observer, considering the changes in the refractive index along the way. Therefore, while Snell's Law quantifies the angle of refraction at the boundary between two mediums, the Fermat principle explains why light follows this particular path by minimizing the time it takes to travel between two points. In media with graded indices, where the refractive index changes spatially, the Fermat principle affects the path of light rays, resulting in curved trajectories explained by Snell's Law.\\

\begin{figure}[h!]
\centering 
\includegraphics[height=0.3\textheight]{{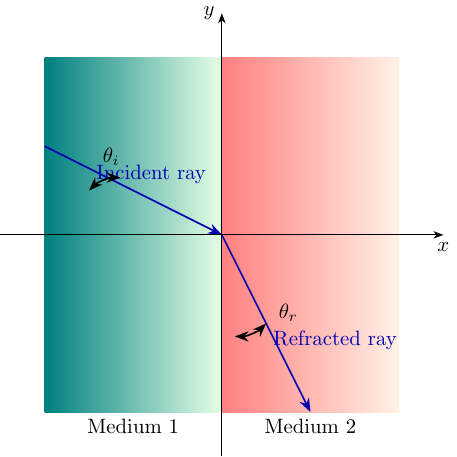}}
\caption{An illustration of a Snell's law.}
\label{fig1:Aa}
\end{figure}

Figure \ref{fig1:Aa} illustrates the fundamental concept of Snell's law \cite{sun2009chebyshev}. In graded-index media, the light path bends or refracts by Snell's law. This problem leads to a more complex radiative transfer solution than uniform-index media. Typically, in such instances, the refractive index varies from the center of the material to its surface.

\noindent
\subsection{The model}
\subsubsection*{Radiation transport equation in graded index media}\label{ab}
The radiation transport equation for graded-index media, which describes the distribution of radiative intensity \( I(s, \boldsymbol{\Omega}, t) \) at position \( s \) and time \( t \) in the direction \( \Omega(\theta, \phi) \), is given as \cite{ymeli2023generalized}.

\begin{equation}\label{eqn:A}
\frac{n}{c_{0}}\dfrac{\partial}{\partial t}I(t,s,\boldsymbol \Omega)+ (k_{e}+ \boldsymbol\Omega.\nabla)I(t,s,\boldsymbol\Omega)+ \frac{1}{n \sin \theta}.I_{\theta} + \frac{1}{n \sin \theta}.I_{\phi} = \mathscr{S}(t,s,\boldsymbol\Omega),
\end{equation}
Where $I_{\theta} =\frac{\partial}{\partial \theta}\left\lbrace I(t,s,\boldsymbol\Omega)(\boldsymbol\Omega \cos \theta - k).\nabla n \right\rbrace $ and $ I_{\phi} = \frac{\partial}{\partial \phi}\left\lbrace (s_{1}.\nabla n)I(t, s,\boldsymbol\Omega)\right\rbrace $. Here, $n$, $s$, and $k_e$  are the refractive index, position vector, and extinction coefficient. $\boldsymbol\Omega = i\sin \theta cos\phi + j \sin \theta \sin \phi + k \cos \theta $ is the direction vector, and $s_{1} = k \times \frac{\boldsymbol\Omega}{\mid k \times \boldsymbol\Omega\mid} =-i\sin \phi +j\cos \phi$, where \( i \), \( j \), and \( k \)   signify the unit vectors in the coordinate system. The sum of coefficients \( k_a \)(absorption) and \( k_s \)(scattering) represent  extinctions coefficients. The sum of $k_a$ and $k_s$ equals the extinction coefficients $k_e$. The single scattering albedo can be expressed as \( \boldsymbol{\omega} = \frac{k_s}{k_e} \). The source term \( \mathscr{S}(s, \boldsymbol{\Omega}, t) \) is provided as \cite{ymeli2023generalized}, \cite{wei2019reverse}
\[ \mathscr{S}(t,s,\boldsymbol\Omega) = n^{2}k_{a}I_{b}(T_{g}) + \frac{k_{s}}{4\pi}\int\limits_{4\pi} \Phi(\Omega' \rightarrow  \Omega)I(t, s, \Omega')d\Omega', \]  
%\[ s(s,\boldsymbol\Omega,t) = n^{2}k_{a}I_{b}(T_{g}) + \frac{k_{s}}{4\pi}\int_{4\pi} \Phi(\Omega' \rightarrow  \Omega)I(s, \Omega', t)d\Omega'  \]  
where $T_g$ denote medium temperature. Scattering phase function \(\Phi(\Omega' \rightarrow \Omega)\) describes the energy redistribution from the incoming direction \(\Omega'(\theta', \phi')\) to the outgoing direction \(\Omega(\theta, \phi)\).\\

Boundary conditions of the model are \(\boldsymbol{\Omega} \cdot \hat{\boldsymbol{n}}_{\boldsymbol{\omega}} < 0\) and \(\boldsymbol{\Omega}_{rs} \cdot \hat{\boldsymbol{n}}_{\boldsymbol{\omega}} < 0\). In this context, \(\hat{\boldsymbol{n}}_{\boldsymbol{\omega}}\) is the normal vector that points outward from the surface. The direction \(\boldsymbol{\Omega}_{rs} = \boldsymbol{\Omega} - 2(\boldsymbol{\Omega} \cdot \hat{\boldsymbol{n}}_{\boldsymbol{\omega}}) \hat{\boldsymbol{n}}_{\boldsymbol{\omega}}\) represents the incident direction for specular reflection. The boundary intensity for $\boldsymbol\Omega_{rs}.\hat{\boldsymbol{n}}_{\boldsymbol{\omega}}<0$ is 
\[ I(t, s_{\boldsymbol{\omega}}, \boldsymbol\Omega) =  Q_{ext} + \rho_{\boldsymbol{\omega}}^{spec} I(s_{t,\boldsymbol{\omega}}, \boldsymbol\Omega_{s})+ n^{2}_{\boldsymbol{\omega}}\epsilon_{\boldsymbol{\omega}}I_{b\boldsymbol{\omega}}+ \frac{\rho_{\boldsymbol{\omega}}^{d}}{\pi}\int\limits_{\boldsymbol\Omega'.\hat{\boldsymbol{n}}_{\boldsymbol{\omega}}>0} \vert \boldsymbol\Omega'.\hat{\boldsymbol{n}}_{\boldsymbol{\omega}}\vert I(t, s_{\boldsymbol{\omega}}, \Omega')d\Omega',  \]
Where \(I_{b\boldsymbol{\omega}}\) represents the black-body intensity at the boundary, while \(\rho^{spec}\) and \(\rho^{d}\) denote the specular and diffuse reflectivities of the boundary, respectively, with emissivity \(\epsilon_{\boldsymbol{\omega}}\). The external driving force \(Q_{ext}\) is incoming from the external side of the boundary at a direction \(\boldsymbol{\Omega}_{0}\) as described by Snell's law. From the RTE, the radiative flux vector \(q(t,s)\) and the incident radiation \(G(t,s)\) at any point \(M(x,y)\) can be computed as follows:

\[
q(t,s) = \int\limits_{4\pi} \boldsymbol{\Omega} I(t,s, \boldsymbol{\Omega}) d\Omega, \quad G(t,s) = \int\limits_{4\pi} I(t,s, \boldsymbol{\Omega}) d\Omega.
\]
Laser irradiation's collimated intensity \(I_{c}\) experiences attenuation as it travels through the medium. This problem can be solved analytically within the medium, with the boundary condition \(I_{c}(t,s_{\omega}, \boldsymbol{\Omega}) = Q_{ext}\). The reduction in collimated intensity \(I_{c}\) in the medium results in the generation of diffuse intensity  \(I_{d}(t,s,\boldsymbol{\Omega})\). Consequently, the total intensity \(I(t,s, \boldsymbol{\Omega})\) comprises both collimated and diffuse components, expressed as:
\[
I(t,s, \boldsymbol{\Omega}) = I_{c}(t, s, \boldsymbol{\Omega}) + I_{d}(t,s,\boldsymbol{\Omega}).
\]

This equation represents the transport equation for graded-index media and is a fundamental framework for modeling radiation transport in such materials.
%\subsection*{Energy redistribution and underlying domain }
\subsection{The underlying model and domain}
Here, $\boldsymbol{\Omega} \in S,  S \subseteq \mathbb{S}^{d-1}$(sphere), $ t \in\left[ 0, T\right] $, position variable $ s \in \textit{D} \subset \mathbb{R}^{d}$, $ \textit{D}_{T} = [0, T] \times \textit{D}$, $I:[0,T]\times \textit{D} \times S  \rightarrow \mathbb{R} $, $ k_{e} = k(s) : \textit{D} \rightarrow  \mathbb{R}_{+}$,  $ k_{a} = k(s) : \textit{D}  \rightarrow  \mathbb{R}_{+}$,  $ k_{s} = k(s) : \textit{D}  \rightarrow  \mathbb{R}_{+}$, $ \Phi : S \times S \rightarrow \mathbb{R}.$\\

Here, the above-defined model will be applied to two approaches: the data-driven approach, also known as the forward problem, and the data-driven discovery approach, commonly referred to as data assimilation or the inverse problem.
\subsubsection{The underlying model for forward problem}
Here, we will define the radiation transport equation with supplemented initial conditions(I.C.) and boundary conditions(B.C.).
The defined partial integrodifferential supplied the I.C. and B.C., which are the following:
\subsection*{I.C.}
\begin{equation}
I(0, s,\boldsymbol\Omega)=I_{0}(s,\boldsymbol\Omega) ~~  \text{where}~~ (s,\boldsymbol\Omega) \in \textit{D} \times  S,  
\end{equation} And $I_{0}:\textit{D}\times S \rightarrow \mathbb{R},$
\subsection*{B.C.}
\begin{equation}
\beta=\left\lbrace  (t,s,\boldsymbol\Omega) \in [0,T] \times \partial{\textit{D}} \times  S: \boldsymbol\Omega.\hat{\boldsymbol{n}}_{\boldsymbol{\omega} }< 0 \right\rbrace. 
\end{equation}
We can wrote $I(t,s,\boldsymbol\Omega)=I_{b}(t,s,\boldsymbol\Omega)$ for some B.C, we can wrote  $I_{b}:\beta \rightarrow \mathbb{R}.$

\subsubsection{The underlying model for inverse problem}
The underlying equation with solution \( I \)  is considered within the subdomain \( \textit{D}'_T \times S  \). This assumption holds that the operator \( \mathscr{L} \) applied to \( I \) in this region equals a given data \( g \). Mathematically, it can be denoted as:
\[
\mathscr{L}(I) = g, \quad \forall (s, \boldsymbol{\Omega}) \in \textit{D}'_T \times S,
\]
Where $\textit{D}' \subset \textit{D}$, $\textit{D}_{T}' = [0,T] \times \textit{D}'$ and g is a source term.

\subsection{Quadrature rules}

Following the approaches outlined in \cite{mishra2021physics, mishra2022estimates, mishra2023estimates}, let \( \boldsymbol D \) represent a domain and \( \varrho \) be an integrable function defined as \( \varrho: \boldsymbol D \rightarrow \mathbb{R} \). Consider the space-time domain \( \textit{D}_{T} = [0,T] \times \textit{D} \subset \mathbb{R}^{d} \), where \( \bar{d} = 2d + 1 \geq 1 \).
We define a mapping \( \varrho: \boldsymbol D \rightarrow \mathbb{R} \), where
\begin{equation}
\varrho = \int_{ \boldsymbol D } \varrho(z) \, dz,
\end{equation}
with \( dz \) representing the \( \bar{d} \)-dimensional Lebesgue measure. To approximate this integral, we utilize quadrature points \( z_{i} \in \boldsymbol D \) for \( 1 \leq i \leq N \), along with their corresponding weights \( w_{i} \) The quadrature approximation of \( \varrho \) is expressed as:
\begin{equation}
\varrho_{N} = \sum_{i=1}^{N} w_{i} \varrho(z_{i}),
\end{equation}
Here, \( z_{i} \) represents the quadrature points. For cases where \( \bar{d} \leq 4 \), standard composite Gauss quadrature rules can be applied using an underlying grid. The selection of quadrature points and weights depends on the order of the quadrature rule, as outlined in \cite{stoer2002r}. However, Gauss quadrature becomes inefficient for higher-dimensional domains.
In cases of moderate dimensionality (\( 4 \leq \bar{d} \leq 20 \)), low-discrepancy sequences, such as Sobol and Halton sequences, are effective for selecting quadrature points, assuming that the function \( \varrho \) has a bounded Hardy-Krause variation \cite{caflisch1998monte}.
In cases of high dimensionality (when \( d \gg 20 \)), Monte Carlo quadrature is the favored approach, exacting randomly chosen quadrature points that are independent and identically distributed(i.i.d.) \cite{caflisch1998monte}.
Let \( \boldsymbol{S} \) represent a training dataset, and define the space-time domain as \( \textit{D}_{T} = [0,T] \times \textit{D} \).
 The selection of training set \( \boldsymbol{S} \subseteq [0,T] \times \textit{D} \) will be depend on appropriate quadrature points. For the PINN algorithms, we choose random points \( z_{i}^{a} = (\boldsymbol{\Omega}_{i}^{\boldsymbol{S}}) \) for \( 1 \leq i \leq N_{\boldsymbol{S}} \), where \( w_{i}^{\boldsymbol{S}} \) represents the corresponding quadrature weights, and \( a \geq 1 \). The Quasi-Monte Carlo (QMC) quadrature method remains unaffected by the curse of dimensionality. In situations where the geometry of the domain is particularly intricate, random points can be chosen as training samples, which, i.i.d. according to the underlying uniform distribution.

\subsection{Training points}\label{Training}
\begin{figure}[htbp]
\centering
\includegraphics[height=0.3\textheight]{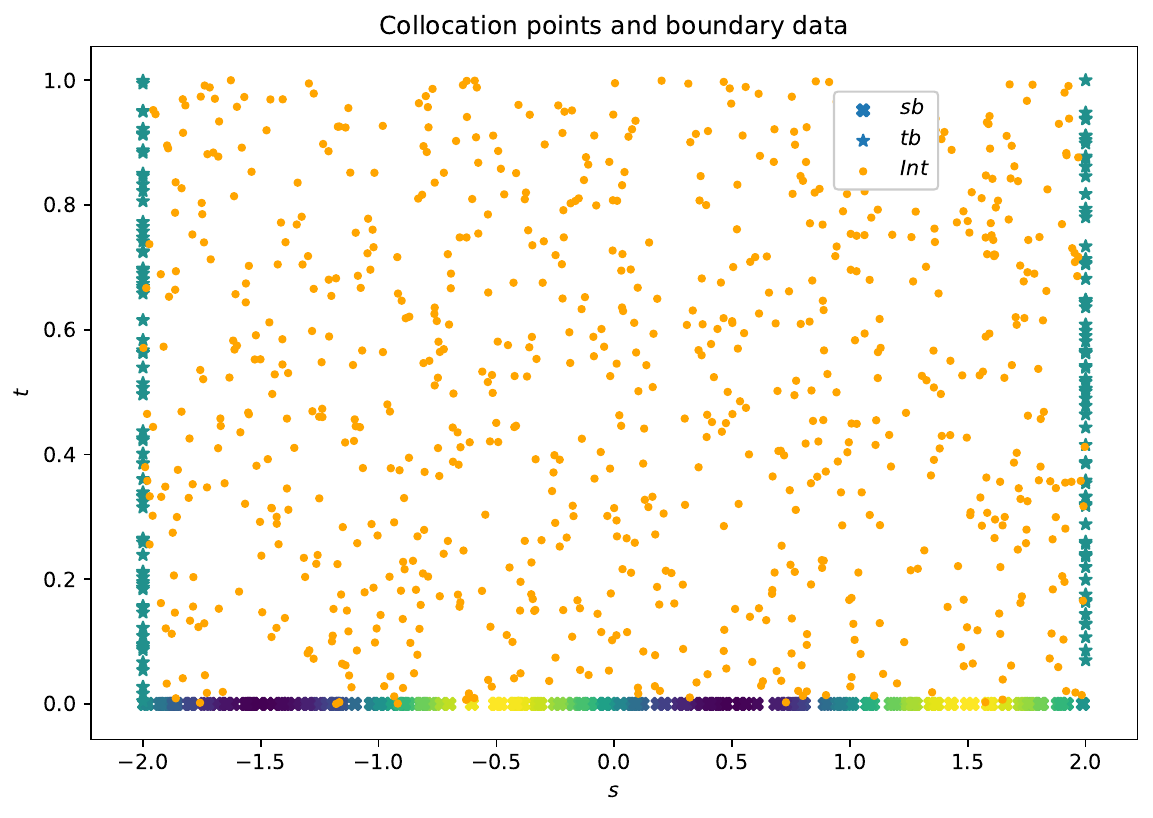}
\caption{Training points (forward problem): A depiction of the training set \( \boldsymbol{S} \) with randomly selected training points. Yellow dots represent interior points, while green and blue dots indicate the temporal and spatial boundary points.}
\label{figt}
\end{figure}

\begin{figure}[htbp]
\centering
\includegraphics[height=0.40\textheight]{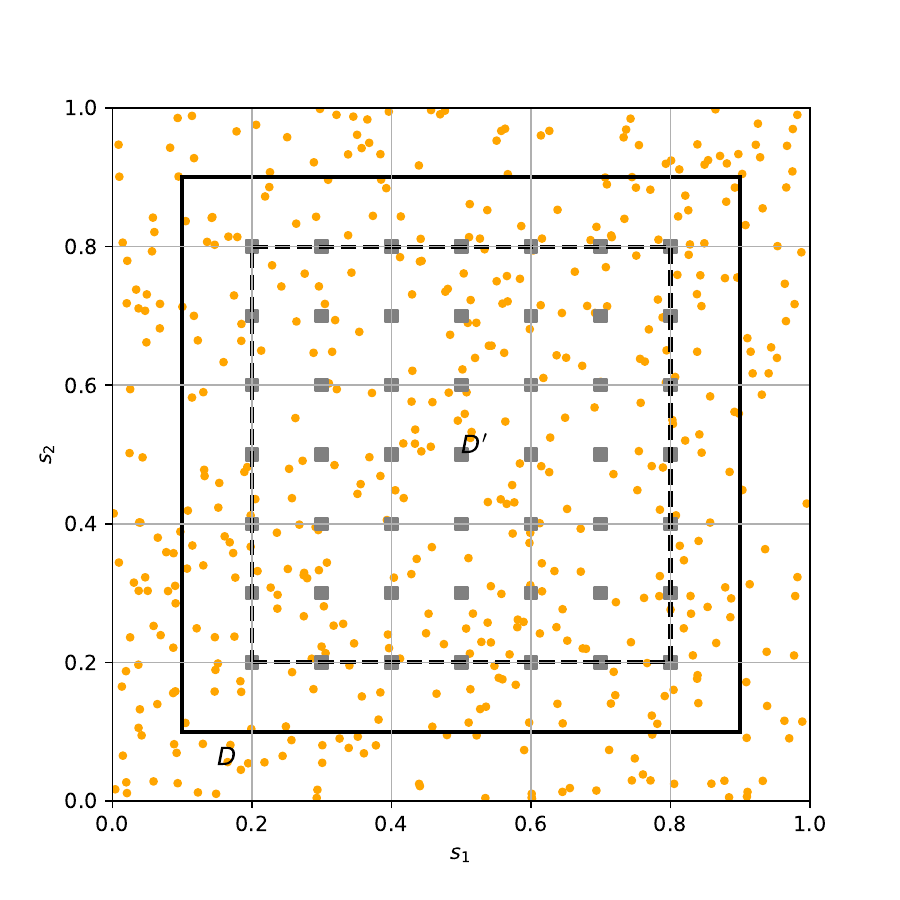}
\caption{Training points (inverse problem): A representation of the training set \( \boldsymbol{S} \) with randomly selected training points. Yellow dots indicate interior points, and grey dots represent Sobol points.}
\label{figtt}
\end{figure}

Physics informed neural networks require four types of training points as described in \cite{mishra2021physics, mishra2022estimates}: interior points \(\zeta_{\text{int}}\), temporal boundary points \(\zeta_{\text{tb}}\), spatial boundary points \(\zeta_{\text{sb}}\), and data points  \(\zeta_{\boldsymbol{d}}\). Figs.\ref{figt} and \ref{figtt} illustrate the training points used in forward and inverse problems (steady state). In Fig. \ref{figtt} $s_{1}$ and $s_{2}$ represent spatial position respectively. 

\subsubsection{Interior training points}
The interior training points are denoted by $\zeta_{\text{int}} = \left\lbrace z_{j}^{\text{int}} \right\rbrace $ for $ 1 \leq j \leq N_{\text{int}} $, where $ z_{j}^{\text{int}} = \left( t_{j}^{\text{int}}, s_{j}^{\text{int}}, \boldsymbol{\Omega}_{j}^{\text{int}}\right) $. Here, $t_{j}^{\text{int}} \in \left[ 0, T \right] $, $s_{j}^{\text{int}} \in D $, and $\boldsymbol{\Omega}_{j}^{\text{int}} \in S $, for all \( j \). These points correspond to quadrature points with weights $w_{j}^{\text{int}} $ based on a suitable quadrature rule. In domains \(\textit{D}\) that are logically rectangular, one can either use Sobol points or randomly select points to create the training set.

\subsubsection{Temporal boundary training points}
The temporal boundary points are represented as \( \zeta_{\text{tb}} = \left\lbrace z_{j}^{\text{tb}} \right\rbrace \), for \( 1 \leq j \leq N_{\text{tb}} \), with \( z_{j}^{\text{tb}} = \left( s_{j}^{\text{tb}}, \boldsymbol{\Omega}_{j}^{\text{tb}} \right) \). Here, \( s_{j}^{\text{tb}} \in D \), and \( \boldsymbol{\Omega}_{j}^{\text{tb}} \in S \), \(\forall\) \( j \). The designated points function as quadrature points within an appropriate quadrature rule, accompanied by weights denoted as \( w_{j}^{\text{tb}} \). For logically rectangular domains \( \textit{D} \), Sobol points can be chosen, or alternatively, random points can be employed to construct the training dataset, similar to the above method.

\subsubsection{Spatial boundary training points}
The spatial boundary points are denoted as \( \zeta_{\text{sb}} = \left\lbrace z_{j}^{\text{sb}} \right\rbrace \), for \( 1 \leq j \leq N_{\text{sb}} \), where \( z_{j}^{\text{sb}} = \left( t_{j}^{\text{tb}}, s_{j}^{\text{tb}}, \boldsymbol{\Omega}_{j}^{\text{tb}} \right) \). In this case, \( t_{j}^{\text{tb}} \in \left[ 0, T \right] \), \( s_{j}^{\text{tb}} \in \partial \textit{D} \), and \( \boldsymbol{\Omega}_{j}^{\text{tb}} \in S \). For logically rectangular domains \( \textit{D} \), Sobol points can be chosen, or alternatively, random points can be employed to construct the training dataset, similar to the above method.
\subsubsection{Data training points}
The data training set is defined as \( \zeta_{\boldsymbol{d}} = \left\lbrace z_{j}^{\boldsymbol{d}} \right\rbrace \) for \( 1 \leq j \leq N_{\boldsymbol{d}} \), where \( z_{j}^{\boldsymbol{d}} \in \textit{D}' \subset \textit{D} \).

\begin{figure}[htbp]
\centering
\includegraphics[height=0.30\textheight]{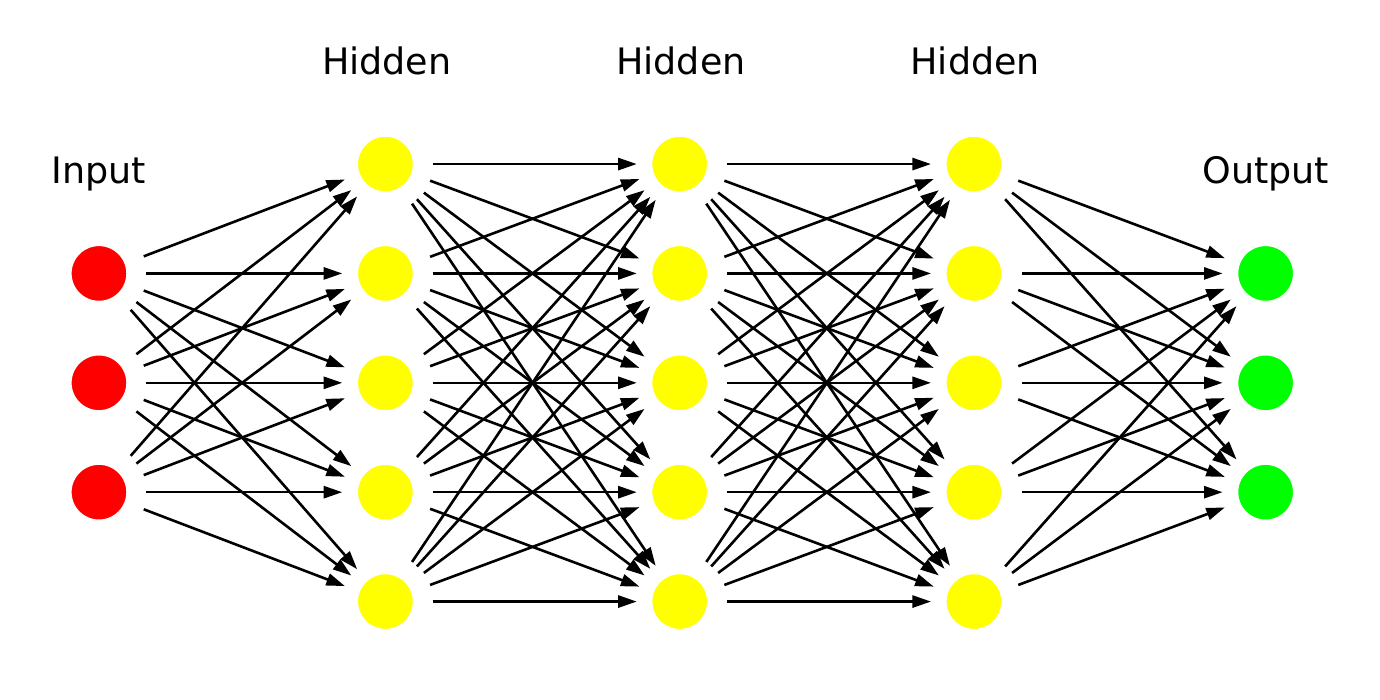}
\caption{In this diagram, input layer neurons are depicted in red, hidden layer neurons in yellow, and output layer neurons in green.}
\label{figttt}
\end{figure}
\subsection{Neural networks}
The PINN functions as a feed-forward neural network, as illustrated in Fig.~\ref{figttt}. A neural network without an activation function behaves as a simple multiple regression model. However, the activation function adds non-linearity, allowing the network to learn and execute more complex tasks. The sigmoid, hyperbolic tangent (tanh), and ReLU functions are some examples of activation functions \cite{goodfellow2016deep}.

The input to the network is \( y = (t, s, \boldsymbol{\Omega}) \in \boldsymbol{D} = [0,T] \times \textit{D} \times S\). The neural network can be expressed as an affine map:
\begin{equation}\label{eqn:N}
I_{\Theta}(y) = C_K \circ  \sigma \circ C_{K-1} \circ \dots \circ \dots    \sigma \circ C_1 (y).
\end{equation}
Where \( \circ \) and  \( \sigma \) are a composition of functions and activation functions, respectively.
For any layer \( k \) such that \( 1 \leq k \leq K \), the transformation at the \( k \)-th layer is defined as follows:
\begin{equation}
C_k z_k = W_k z_k + b_k \quad \text{where} \quad W_k \in \mathbb{R}^{d_{k+1} \times d_k}, \, z_k \in \mathbb{R}^{d_k}, \quad \text{and} \quad b_k \in \mathbb{R}^{d_{k+1}}.
\end{equation}
To ensure consistency, we establish \( d_1 = \bar{d} = 2d + 1 \), where \( d \) represents the spatial dimension, and we define \( d_K = 1 \) for the output layer. In the context of machine learning, this neural network consists of an input layer, an output layer, and \( K-1 \) hidden layers, with the condition that \( 1 < K < \mathbb{N} \).
Each hidden layer \( k \), which consists of \( d_k \) neurons, takes an input vector \( z_k \in \mathbb{R}^{d_k} \). Initially, this input vector is transformed using the linear mapping \( C_k \), after which it is processed by a non-linear activation function represented as \( \sigma \). The overall count of neurons in the network can be expressed as \( 2d + 2 + \sum_{k=2}^{K-1} d_k \).
We define the parameter set for the network, which includes weights and biases, as \( \Theta = \left\lbrace W_k, b_k \right\rbrace \). Furthermore, we denote the set of weights as \( \Theta_w = \left\lbrace W_k \right\rbrace \) for all \( 1 \leq k \leq K \)  \cite{mishra2021physics, mishra2022estimates}.
The parameters \( \Theta \) lie in the space \( \Theta' \subset \mathbb{R}^P \), and \( P \) denote total number of parameters:
\begin{equation}
P = \sum_{k=1}^{K-1} \left( d_k + 1 \right) d_{k+1}.
\end{equation}
\subsection{Residuals}
This section outlines the residuals associated with the divided training sets, including interior, temporal, data (for inverse), and spatial training points. The focus is on minimizing these residuals. Optimization techniques will involve stochastic gradient descent methods, such as ADAM for first-order optimization, and higher-order approaches like variations of the BFGS algorithm. The PINN \( I_{\Theta} \) depends on tuning parameters \( \Theta \in \Theta' \), which represent the weights and biases within the network. In a typical deep learning setup, the network is trained by optimizing these parameters \( \Theta \) to ensure that the neural network approximation \( I_{\Theta} \) accurately aligns with the exact solution \( I \). The interior residual can be written as:
\begin{equation}
\mathrm{R}_{\text{int},\Theta} = \mathrm{R}_{\text{int},\Theta}(t,s,\boldsymbol\Omega), \quad \forall (t,s,\boldsymbol\Omega) \in [0,T] \times \textit{D} \times S,
\end{equation}
We can express the interior residual \( \mathrm{R}_{\text{int}, \Theta} \) as:
\begin{equation}\label{eqn:Ra}
\begin{split}
\mathrm{R}_{\text{int}, \Theta} &= \frac{n}{c_{0}} \frac{\partial}{\partial t}I_{\Theta} + (k_{e} + \boldsymbol{\Omega} \cdot \nabla)I_{\Theta} + \frac{1}{n \sin \theta} \frac{\partial}{\partial \theta} \left\lbrace I_{\Theta} (\boldsymbol{\Omega} \cos \theta - k) \cdot \nabla n \right\rbrace \\
&\quad + \frac{1}{n \sin \theta} \frac{\partial}{\partial \phi} \left\lbrace (s_{1} \cdot \nabla n) I_{\Theta} \right\rbrace - n^{2} k_{a} I_{b \Theta}(T_{g})- \frac{k_{s}}{4\pi} \sum_{i=1}^{N_{\boldsymbol{S}}} w_{i}^{\boldsymbol{S}} \Phi(\boldsymbol{\Omega}, \boldsymbol{\Omega}_{i}^{\boldsymbol{S}}) I_{\Theta},
\end{split}
\end{equation}
In this equation, \( (\boldsymbol{\Omega}_{i}^{\boldsymbol{S}}) \) represent the Gauss-Legendre quadrature points, while \( w_{i} \) signify the corresponding quadrature weights of order \(\boldsymbol{S} \). The residuals for the initial, boundary, and data points are defined as follows:
\begin{equation}\label{eqn:Rab}
\begin{aligned}
\mathrm{R}_{\text{tb}} &= \mathrm{R}_{\text{tb}, \Theta} = I_{\Theta} - I_{0}, \quad \forall (s, \boldsymbol{\Omega}) \in \textit{D} \times S, \\
\mathrm{R}_{\text{sb}} &= \mathrm{R}_{\text{sb}, \Theta} = I_{\Theta} - I_{b}, \quad \forall (t, s, \boldsymbol{\Omega}) \in \beta.
\end{aligned}
\end{equation}
And
\begin{equation}\label{eqn:Residual_data}
\begin{aligned}
\mathrm{R}_{ {\boldsymbol{d}}} &= \mathscr{L}(I_{\Theta}) - g, \quad \forall (s, \boldsymbol{\Omega}) \in \textit{D}'_{T} \times S.  
\end{aligned}
\end{equation}
We aim to find the optimal set of tuning parameters \( \Theta \in \Theta' \) that minimizes the residual specified in the forward problem equation,
\begin{equation}\label{eqn:R}
  \Theta^{\ast} \in \Theta' : \Theta^{\ast} = \arg \min_{\Theta \in \Theta'} \left( \Vert \mathrm{R}_{\text{int}, \Theta} \Vert^{2}_{L^{2}([0,T] \times \textit{D} \times S)} + \Vert \mathrm{R}_{\text{sb}, \Theta} \Vert^{2}_{L^{2}(\beta)} + \Vert \mathrm{R}_{\text{tb}, \Theta} \Vert^{2}_{L^{2}(\textit{D} \times S )} \right).
\end{equation}
For the inverse problem, we add the term corresponding to the data residual \( \mathrm{R}_{\boldsymbol{d}} \) to Eq.\eqref{eqn:R}. This results in the following minimization problem:
\begin{equation}\label{eqn:Rb}
  \Theta^{\ast} \in \Theta' : \Theta^{\ast} = \arg \min_{\Theta \in \Theta'} 
  \left( 
  \Vert \mathrm{R}_{\text{int}, \Theta} \Vert^{2}_{L^{2}( \textit{D}_{T} \times S)} 
  + \Vert \mathrm{R}_{\text{sb}, \Theta} \Vert^{2}_{L^{2}(\beta)}  
  + \Vert \mathrm{R}_{\boldsymbol{d}, \Theta} \Vert^{2}_{L^{2}( \textit{D}'_{T} \times S)}  
  \right).
\end{equation}
The integrals in Eqs.\eqref{eqn:R} and \eqref{eqn:Rb} cannot be calculated exactly due to the use of the \( L^{2} \) norm, necessitating approximation through an appropriate quadrature method.

\subsection{Loss functions and optimization}\label{loss}
We approximate the above integral with the following loss functions for forward and inverse problems, respectively
\begin{equation}
\mathscr{J}_{1}(\Theta)=\sum\limits_{j=1}^{N_{sb}}w_{j}^{sb}\vert  \mathrm{R}_{\text{sb}, \Theta} (z_{j}^{sb}) \vert^{2}+ \sum\limits_{j=1}^{N_{tb}}w_{j}^{tb}\vert  \mathrm{R}_{\text{tb}, \Theta} (z_{j}^{tb}) \vert^{2}+\lambda \sum\limits_{j=1}^{N_{int}}w_{j}^{int}\vert  \mathrm{R}_{\text{int}, \Theta} (z_{j}^{int}) \vert^{2} \label{eqn:La},
\end{equation}
\begin{equation}
\mathscr{J}_{2}(\Theta)=\sum\limits_{j=1}^{N_{d}}w_{j}^{d}\vert  \mathrm{R}_{\boldsymbol{d}, \Theta} (z_{j}^{d}) \vert^{2}+\sum\limits_{j=1}^{N_{sb}}w_{j}^{sb}\vert  \mathrm{R}_{\text{sb}, \Theta} (z_{j}^{sb}) \vert^{2}+ \lambda \sum\limits_{j=1}^{N_{int}}w_{j}^{int}\vert  \mathrm{R}_{\text{int}, \Theta} (z_{j}^{int}) \vert^{2} \label{eqn:Lb},
\end{equation}
Regularize the minimization problems for the loss function, i.e
\begin{equation}\label{eqn:L}
\Theta^{\ast} = \arg \min_{\Theta \in \Theta'}(\mathscr{J}_{i}(\Theta)+ \lambda_{reg}\mathscr{J}_{reg}(\Theta)),
\end{equation}
Where \( i = 1, 2 \). In machine learning, it is common to incorporate a regularization term to mitigate overfitting. A widely used form of the regularization function is \( \mathscr{J}_{\text{reg}}(\Theta) = \Vert \Theta \Vert_q^q \), where \( q \) is typically 1 (for \( L^{1} \) regularization) or 2 (for \( L^{2} \) regularization).
The parameter \( \lambda_{\text{reg}} \) controls the balance between the regularization term and the actual loss function \( \mathscr{J} \), with \( 0 \leq \lambda_{\text{reg}} \ll 1 \).
Stochastic gradient descent methods, including ADAM, will be applied due to their popularity in first-order optimization. Additionally, advanced optimization techniques, including various forms of the BFGS algorithm, may be employed.
Our aim is to identity the optimal solution \( I^{\ast} = I_{\Theta^{\ast}} \) using the training sets. We begin with an initial value \( \bar{\Theta} \in \Theta' \) and calculate the network output \( I_{\bar{\Theta}} \), PDE residual, boundary residual, loss function, and its gradients. Ultimately, the optimal solution is \( I^{\ast} = I_{\Theta^{\ast}} \), which is determined by the PINN.\\

We approximate local minimum in \eqref{eqn:L} as \( \Theta^{\ast} \). The resulting DNNs \( I^{\ast} = I_{\Theta^{\ast}} \)  will  solution \( I \) of \eqref{eqn:A}.
The Table\ref{table_1} contains a hyperparameter of numerical experiments.
We summarize the PINN algorithms for approximating RTE in a graded index medium. The algorithms are  described in \cite{mishra2021physics}, \cite{mishra2022estimates}, \cite{bai2021physics}, and \cite{mishra2023estimates}. Below, Algorithm \ref{alg1} is presented for forward problems, while Algorithm \ref{alg2} addresses inverse problems:
\begin{table}[!ht]
\centering
\caption{The configurations of hyperparameters and the frequency of retraining utilized in ensemble training for physics-informed neural networks (PINN).}
\begin{tabular}{||c c c c c c||}
\hline
Examples & $K-1$ & $\bar{d}$ & $\lambda$ & $\lambda_{\text{reg}}$ & $n_{\Theta}$ \\ [0.5ex] 
\hline
Example 1a,b,c & 4, 8 & 20, 24 & 0.1, 1, 10 & 0 & 4 \\ 
\hline
~~~Example 2a,b,c,d & 4, 8 & 24, 28 & 0.1, 1, 10 & 0 & 10,10,4,4 \\ 
\hline
Example 3a,b,c & 4, 8 & 24, 28 & 0.1, 1, 10 & 0 & 4 \\ 
\hline
Example 4a,~b & 4, 8 & 20, 24 & 0.1, 1, 10 & 0 & 4 \\ 
\hline
Example 5a,b,c & 4, 8 & 20, 24 & 0.1, 1, 10 & 0 & 4 \\ 
\hline
Example 6a,~b & 4, 8 & 20, 24 & 0.1, 1, 10 & 0 & 4,5 \\ 
\hline
\end{tabular}
\label{table_1}
\end{table}

\begin{algorithm} 
\label{alg1} {\bf Algorithm for developing a PINN to estimate radiative intensity in forward problems.}
\begin{itemize}
\item [{\bf Inputs}:] Underlying domain, data, and coefficients for the RTE with graded index Eq.(\ref{eqn:A}); quadrature points and weights for the underlying quadrature rules; non-convex gradient-based optimization algorithms.
\item [{\bf Aim}:] To approximate the solution of the model, using a PINN \(I^{\ast} = I_{\Theta^{\ast}}\).
\item [{\bf Step $1$}:] Select the training sets as outlined in section \ref{Training}.
\item [{\bf Step $2$}:] Initialize with a weight vector \( \bar{\Theta} \in \Theta' \) and calculate: neural network \( I_{\bar{\Theta}} \) Eq.(\ref{eqn:N}),PDE residual Eq.(\ref{eqn:Ra}), boundary residuals Eq.(\ref{eqn:Rab}), loss function Eq.(\ref{eqn:La}), Eq.(\ref{eqn:L}),
        and gradients for optimization algorithm initiation.
\item [{\bf Step $3$}:] Execute the optimization algorithm until reaching an approximate local minimum \( \Theta^{\ast} \) of Eq.(\ref{eqn:L}). The resulting function \( I^{\ast} = I_{\Theta^{\ast}} \) is the desired PINN for approximating the radiative transfer equation solution \( I \). 
\end{itemize}
\end{algorithm}

\subsection{Estimation on generalization error} 
Let the spatial domain be \(\textit{D} = [0, 1]^{d}\), where \(d\) denote the spatial dimension. This section focuses on obtaining an accurate estimation of the generalization error or so-called total error for the trained neural network, \(I^{\ast} = I_{\Theta^{\ast}}\). This result arises from the application of the PINNs algorithms \ref{alg1} and \ref{alg2}. The error can be expressed as follows:

\begin{equation}\label{eqn:ge}
\mathrm{E}_{G}=\mathrm{E}_{G}(\theta^{\ast})= \left( \int\limits_{\boldsymbol{D}} \vert I(t,s,\boldsymbol{\Omega}) - I^{\ast}(t,s,\boldsymbol{\Omega}) \vert^{2}dX \right)^{\frac{1}{2}}
\end{equation}
Where \( dX = dt \, ds \, d\boldsymbol{\Omega} \) denotes the volume measure on \(\boldsymbol{D}\). This approach outlined in \cite{mishra2021physics}, \cite{bai2021physics}, \cite{mishra2022estimates}, and \cite{mishra2023estimates}. This section provides an estimation of the generalization error, as defined in equation (\ref{eqn:ge}), based on the training error.

\begin{equation}\label{eqn:trans} 
\begin{split}
\mathrm{E}_{T}^{N_{\text{sb}}} &= \left( \sum\limits_{j=1}^{N_{sb}} w_{j}^{sb} \left| \mathrm{R}_{\text{sb}, \Theta^{\ast}} (z_{j}^{sb}) \right|^{2} \right)^{\frac{1}{2}}, 
\mathrm{E}_{T}^{N_{\text{tb}}} = \left( \sum\limits_{j=1}^{N_{tb}} w_{j}^{tb} \left| \mathrm{R}_{\text{tb}, \Theta^{\ast}} (z_{j}^{tb}) \right|^{2} \right)^{\frac{1}{2}}, \\
\mathrm{E}_{T}^{N_{\text{int}}} &= \left( \sum\limits_{j=1}^{N_{int}} w_{j}^{int} \left| \mathrm{R}_{\text{int}, \Theta^{\ast}} (z_{j}^{int}) \right|^{2} \right)^{\frac{1}{2}}, 
\mathrm{E}_{T}^{N_{\boldsymbol{d}}} = \left( \sum\limits_{j=1}^{N_{\boldsymbol{d}}} w_{j}^{\boldsymbol{d}} \left| \mathrm{R}_{\boldsymbol{d}, \Theta^{\ast}} (z_{j}^{\boldsymbol{d}}) \right|^{2} \right)^{\frac{1}{2}}.
\end{split}
\end{equation}
The generalized error is similar to the form presented in \cite{mishra2021physics} and is expressed as:
\begin{equation}
\begin{split}
(\mathrm{E}_{G})^2  & \leq V\left(( \mathrm{E}_{T}^{tb})^{2} + v( \mathrm{E}_{T}^{sb})^{2} 
+ c ( \mathrm{E}_{T}^{int})^{2}\right) \\
& + VV_{2}\left( \dfrac{(\log(N_{tb}))^{2d}}{N_{tb}} + c \dfrac{(\log(N_{sb}))^{2d}}{N_{sb}} 
+ c \dfrac{(\log(N_{int}))^{2d+1}}{N_{int}} + c N_{\boldsymbol{S}}^{-2a} \right).
\end{split}
\end{equation}
Where $V_{2},v,$ and $V$ are constants as defined in  \ref{thm:1}.\\
Consider $\Phi$  is symmetric such that 
\[\Phi(\boldsymbol{\Omega}, \boldsymbol{\Omega}')= \Phi(\boldsymbol{\Omega}', \boldsymbol{\Omega}),\]
And let 
 \[ \Sigma_{g}(\boldsymbol{\Omega})= \int\limits_{4\pi}\Phi(\boldsymbol{\Omega},
 \boldsymbol{\Omega}')d\Omega' \]
Where $S \in \Omega' $ and
$\Sigma_{g}$  is essentially bounded i.e $  \Sigma_{g} \in L^{\infty}(S).$

\begin{algorithm} 
\label{alg2} {\bf Algorithm for developing a PINN to estimate radiative intensity in inverse problems.}
\begin{itemize}
\item [{\bf Inputs}:] Underlying domain, data, and coefficients for the RTE with graded index Eq.(\ref{eqn:A}); appropriate quadrature points and weights for the underlying quadrature rules; and non-convex gradient-based optimization algorithms.
\item [{\bf Aim}:] To approximate the solution \(I\)of Eq.(\ref{eqn:A}) for inverse problems, using a PINN \(I^{\ast} = I_{\Theta^{\ast}}\).
\item [{\bf Step $1$}:] Select the training sets outline in Section~\ref{Training}.
\item [{\bf Step $2$}:] Initialize with a weight vector \( \bar{\Theta} \in \Theta' \) and calculate: neural network \( I_{\bar{\Theta}} \) Eq.(\ref{eqn:N}),PDE residual Eq.(\ref{eqn:Ra}), data residuals Eq.(\ref{eqn:Residual_data}), loss function Eq. (\ref{eqn:Lb}), Eq.(\ref{eqn:L})
        and gradients for optimization algorithm initiation.
\item [{\bf Step $3$}:]Execute the optimization algorithm until reaching an approximate local minimum \( \Theta^{\ast} \) of Eq.(\ref{eqn:L}). The resulting function \( I^{\ast} = I_{\Theta^{\ast}} \) is the desired PINN for approximating the radiative transfer equation solution \( I \). 
\end{itemize}
\end{algorithm}

\subsection{Steady case}
If \( c_{0} \rightarrow \infty \), Eq.\eqref{eqn:A} becomes
\begin{equation}\label{eqn:BB}
(k_{e}+\boldsymbol\Omega.\nabla)I+ \frac{1}{n \sin \theta}.I_{\theta} + \frac{1}{n \sin \theta}.I_{\phi} =  \mathscr{S},
\end{equation}
\subsection*{B.C.}
\begin{equation}
\beta_{o}=\left\lbrace  (s,\boldsymbol\Omega) \in  \partial{\textit{D}} \times  S: \boldsymbol\Omega.\hat{\boldsymbol{n}}_{\boldsymbol{\omega }}< 0 \right\rbrace, 
\end{equation}
We can write $I(s,\boldsymbol\Omega)=I_{b}(s,\boldsymbol\Omega,\xi)$ for some B.C., we can wrote  $I_{b}:\beta_{o} \rightarrow \mathbb{R}.$
Let \(\textit{D} = [0,1]^{d}\), where d denotes spatial dimension. This section aims to derive a precise generalization error estimate for the trained neural network \( I^\ast = I_{\Theta^\ast} \). This network is the result of the PINNs algorithms described in algorithms \ref{alg1} and \ref{alg2}.
\begin{equation}
\mathrm{E}_{G_{steady}}=\mathrm{E}_{G}(\theta^{\ast})= \left( \int\limits_{\textit{D} \times  S} \vert I(s,\boldsymbol{\Omega}) - I^{\ast}(s,\boldsymbol{\Omega}) \vert^{2}dX \right)^{\frac{1}{2}},
\end{equation}
Where $ dX=ds d\boldsymbol{\Omega},  $
\begin{equation}\label{eqn:steady}
\begin{split}
\mathrm{E}_{T}^{N_{\text{sb}}} = \left( \sum\limits_{j=1}^{N_{sb}} w_{j}^{sb} \left| \mathrm{R}_{\text{sb}, \Theta^{\ast}} (z_{j}^{sb}) \right|^{2} \right)^{\frac{1}{2}}, 
\mathrm{E}_{T}^{N_{\text{int}}} &= \left( \sum\limits_{j=1}^{N_{int}} w_{j}^{int} \left| \mathrm{R}_{\text{int}, \Theta^{\ast}} (z_{j}^{int}) \right|^{2} \right)^{\frac{1}{2}}, \\
\mathrm{E}_{T}^{N_{\boldsymbol{d}}} &= \left( \sum\limits_{j=1}^{N_{\boldsymbol{d}}} w_{j}^{\boldsymbol{d}} \left| \mathrm{R}_{\boldsymbol{d}, \Theta^{\ast}} (z_{j}^{\boldsymbol{d}}) \right|^{2} \right)^{\frac{1}{2}}.
\end{split}
\end{equation}
The generalized error is similar to the form presented in \cite{mishra2021physics} and is expressed as:
\begin{align}
(\mathrm{E}_{G_{steady_{f}}})^{2} &\leq V \left( \nu (\mathrm{E}_{T}^{sb})^{2} + \nu (\mathrm{E}_{T}^{int})^{2} \right) \notag \\
&+ V \left( \dfrac{(\log(N_{sb}))^{2d}}{N_{sb}} + \nu \dfrac{(\log(N_{int}))^{2d}}{N_{int}} + \nu N_{\boldsymbol{S}}^{-2a}\right),
\end{align}
Where $v$ and $V$ are constants as defined in \ref{thm:generalization_error_forward_inverse}.\\
The generalization error bound expresses that the approximation error for the underlying problems using a trained PINN will remain small if certain conditions are satisfied:
\begin{remark}
The PINN should be effectively trained, as evidenced by a sufficiently small training error.
Although the training error cannot be controlled \emph{a priori}, it can be computed \emph{a posteriori}. Sufficient training (collocation) points are required to ensure accurate learning.
The quadrature error, which depends on the number of collocation points \(N\) and the quadrature constants, can be reduced by selecting a large enough \(N\). This finding underscores that the generalization error estimate sets an upper bound on the total error, encompassing both training errors (from equations \ref{eqn:trans} and \ref{eqn:steady}) and the number of training data points \(N_{\text{int}}, N_{\text{sb}}, N_{\text{tb}}\), as well as the quadrature points \(N_S\) utilized to approximate the scattering integral (from equation \ref{eqn:A}). Although no \emph{a priori} estimate is available for the training errors, they can be computed after complete training. Therefore, the theorems suggest that the model will provide less relative error if the involved constants remain finite and the PINN is trained adequately. This theory aligns with general machine learning theory, where a well-trained and regularized PINN \(I^*\) ensures stability and bounded generalization error. Here $N_{int}>128, \text{and} ~N_{sb,tb}>64$.
\end{remark}
\begin{remark}
We address the following inverse problem: suppose that the boundary conditions in equations (\ref{eqn:A}) and (\ref{eqn:BB}), which may also include initial conditions, are unknown. This problem is linked to the PDE being an ill-posed problem. However, we assume noiseless measurements of the underlying solution \(I\) are available within a subdomain \(\textit{D}' \subset \textit{D}\). In most cases, PINNs achieve minimal errors with less than one minute of training time. Due to their simplicity and efficiency, PINNs provide an attractive alternative to current data assimilation methods, particularly in high-dimensional problems.
\end{remark}

The error estimate primarily relates the overall generalization error to the training error. It uses the stability of the PDE to establish an upper bound based on the PDE residual, which is influenced by both training and quadrature errors.
As long as the training errors remain independent of the underlying dimensionality, the above estimate implies that the PINNs described in Algorithms \ref{alg1} and \ref{alg2} will not be impacted by the curse of dimensionality.
\section{Numerical experiments} \label{sec:3}
The PINN algorithms \ref{alg1} and \ref{alg2} were implemented using the PyTorch framework \cite{paszke2017automatic}. All numerical experiments were performed on an Apple MacBook with an M3 chip and 24 GB of RAM.
Several key hyperparameters are essential to the PINNs framework, including the quantity of hidden layers \(K - 1\), layer width, selected activation function \( \sigma \), parameter \(\lambda\) in the loss function, regularization parameter \(\lambda_{\text{reg}}\) in the cumulative loss, and the specific gradient descent optimization algorithm. For the activation function \( \sigma \), we select the hyperbolic tangent (\(\tanh\)), which ensures the smoothness properties required by theoretical guarantees for the neural networks are obtained. 
We employ the second-order LBFGS optimizer to improve convergence. For fine-tuning the remaining hyperparameters, we adopt the ensemble training method described in \cite{bai2021physics}, \cite{mishra2021physics}, \cite{mishra2022estimates} and \cite{mishra2023estimates}. This method involves evaluating different values for the number of hidden layers, layer depth, parameter \(\lambda\), and regularization term \(\lambda_{\text{reg}}\), as illustrated in Table \ref{table_1}. Each hyperparameter set is used to retrain the model \(n_\theta\) times in parallel with different random weight initializations. The configuration yielding the lowest training loss is then selected as the best model.
\subsection{Forward problem: data driven}
\subsubsection{Radiation distribution of 1D infinite wall with Gaussian source }\label{sec:1D}
This model is based on a non-scattering medium situated between 1D infinite parallel black walls. The model contains Gaussian source term. We assume a constant extinction coefficient for the medium. The following radiative transfer equation describes this problem \cite{liu2006finite}: 
\begin{equation}\label{eq:1d}
\mu \frac{dI}{dx}+ k_{e}I = \exp\left(-(x-c)^{2}/ \alpha^{2}\right), \quad x, c \in [0,1]
\end{equation}
The boundary conditions are:
\begin{subequations} 
\begin{align}
I(0,\mu) &= k_{e}^{-1}\exp\left(-c^{2}/\alpha^{2}\right), \quad \mu > 0 \label{eq:a} \\
I(1,\mu) &= k_{e}^{-1}\exp\left(-(1-c)^{2}/\alpha^{2}\right), \quad \mu < 0 \label{eq:b}
\end{align}
\end{subequations}
The analytical solution (for \( \mu > 0 \)) is provided as:
\begin{equation}
\begin{split}
I(x,\mu) &= I(0,\mu)\exp\left(- \frac{k_{e}x}{\mu}\right) - \frac{\alpha\sqrt{\pi}}{2\mu} \exp\left\lbrace -\frac{k_{e}}{\mu}\left[ x - \left( \frac{\alpha^{2}k_{e}}{4\mu} +c \right)\right]\right\rbrace \\
 & \times \left[ \operatorname{erf}\left( -\frac{k_{e}\alpha}{2\mu} + \frac{c-x}{\alpha}\right) - \operatorname{erf}\left( -\frac{k_{e}\alpha}{2\mu} + \frac{c}{\alpha}\right)\right].
\end{split}
\end{equation}
\begin{figure}[htbp]
\centering
\includegraphics[height=0.3\textheight]{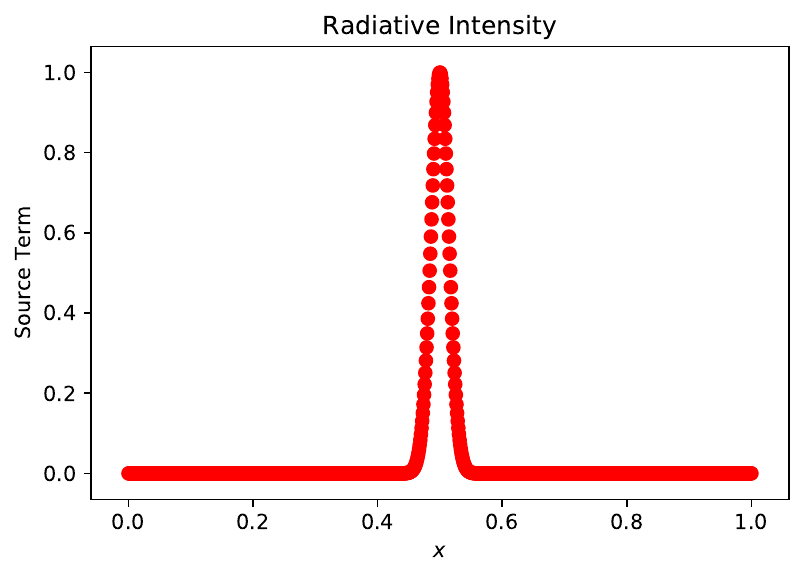}
\caption{Source term of Eq.\ref{eq:1d}.}
\label{source1}
\end{figure}
In this scenario, \( \alpha = 0.02 \), \( c = 0.5 \), and \( \mu = 0.5 \). We compare the radiation intensity distributions calculated for the direction \( \mu = 0.5 \) across media with varying extinction coefficients: \( k_e = 0.1\), \( k_e = 1\), and \( k_e = 10 \, m^{-1}\). Figure \ref{source1} represents the source term of the model. The predicted and exact solutions of the RTE for \( k_e = 0.1\), \( k_e = 1\), and \( k_e = 10 \, m^{-1}\) are presented in Figures \ref{ka_0.1}, \ref{ka_1}, and \ref{ka_10}, respectively. \\
The PINN solution is consistent and closely aligned with the exact solutions for all three scenarios. Table \ref{table_2} demonstrates that the errors remain minimal, further highlighting the PINN's capability to approximate the PDE with low computational cost accurately. 
The authors simulated the model's relative error for parameters \( \mu = 0.5 \), \( \alpha = 0.02 \), and \( c = 0.5 \). For the meshfree method by Zhao et al. \cite{zhao2013second}, the relative error at \( k_{e} = 0.1, 1 \), and \( 10\,m^{-1} \), with meshes ranging from 10 to 400, decreases from \( 21.77 \) to \( 1.95 \times 10^{-3} \), \( 6.35 \) to \( 0.36 \times 10^{-2} \), and \( 2.29 \) to \( 0.49 \times 10^{-2} \), respectively. For the DFEM\cite{feng2018discontinuous}, the relative error at \( k_{e} = 0.1, 1 \), and \( 10\,m^{-1} \), with meshes ranging from 10 to 400, reduces from \( 1.95 \) to \( 3.49 \times 10^{-4} \), \( 2.05 \) to \( 5.09 \times 10^{-4} \), and \( 1.74 \) to \( 0.15 \times 10^{-2} \), respectively. For the generalized lattice Boltzmann method, with meshes between 20 and 200, the relative error at \( k_{e} = 0.1, 1 \), and \( 10\,m^{-1} \) decreases from \( 1.6 \times 10^{-2} \) to \( 8.9 \times 10^{-5} \), \( 7.4 \times 10^{-2} \) to \( 5.3 \times 10^{-4} \), and \( 1.86 \) to \( 0.32 \times 10^{-4} \), respectively. In contrast, the PINN method achieves a relative \( L^{2} \) error of \( 6.4 \times 10^{-4} \) at \( k_{e} = 10\,m^{-1} \).

\begin{figure}[htbp]
\centering
\includegraphics[height=0.3\textheight]{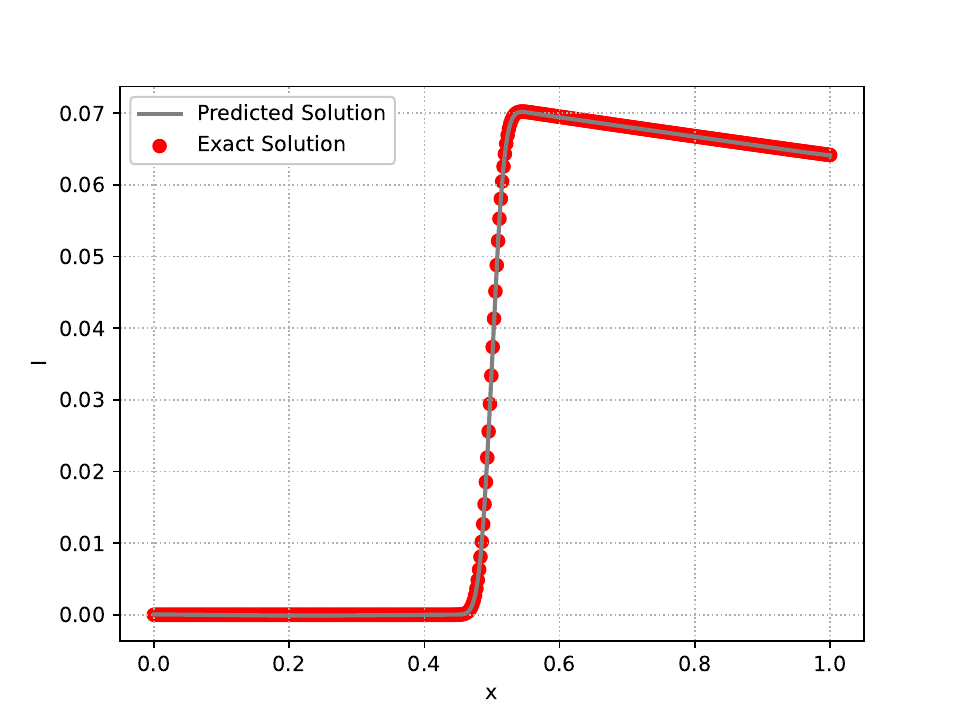}
\caption{Radiation distribution with a Gaussian source at an  extinction coefficients $k_e = 0.1m^{-1}$.}
\label{ka_0.1}
\end{figure}
\begin{figure}[htbp]
\centering
\includegraphics[height=0.3\textheight]{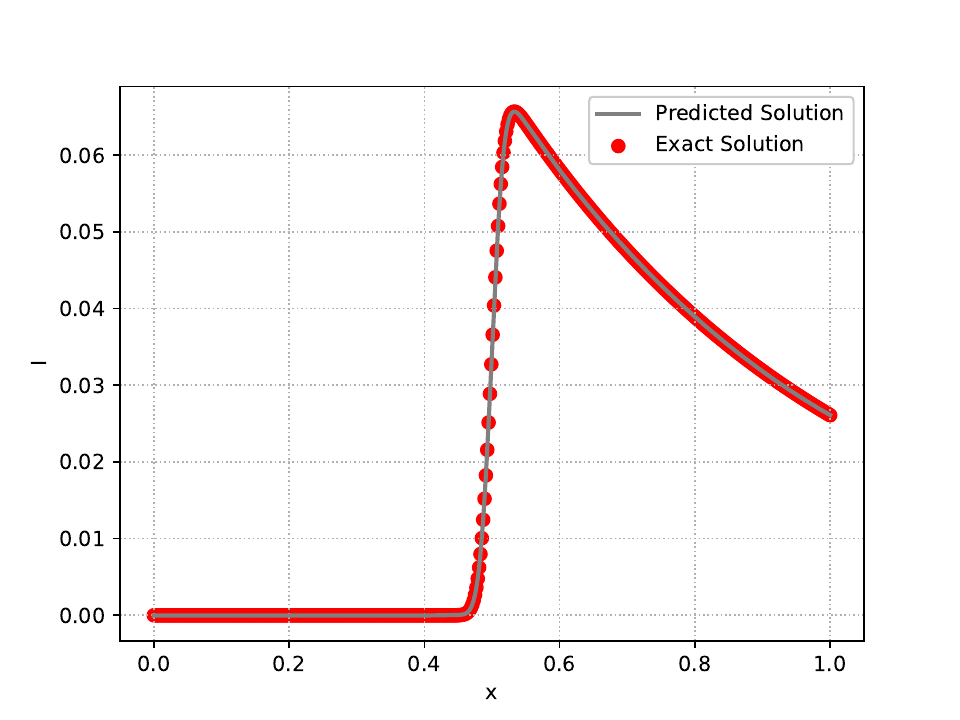}
\caption{Radiation distribution with a Gaussian source at an extinction coefficients $k_e = 1 m^{-1}$.}
\label{ka_1}
\end{figure}
\begin{figure}[htbp]
\centering
\includegraphics[height=0.3\textheight]{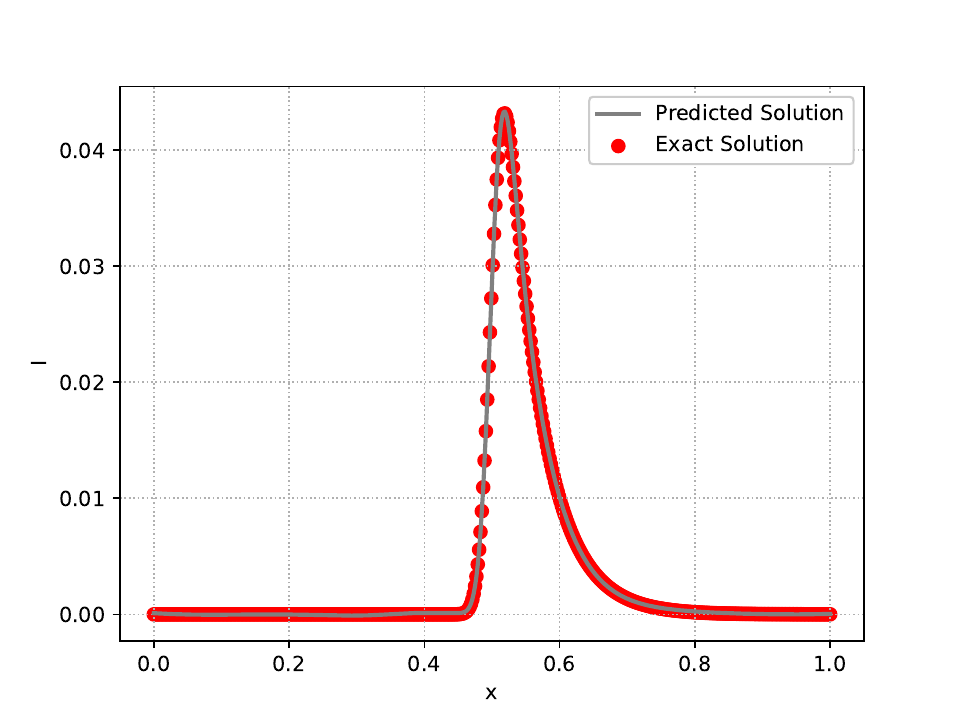}
\caption{Radiation distribution with a Gaussian source at an extinction coefficients $k_e =10 m^{-1}$.}
\label{ka_10}
\end{figure}
\begin{table}[!ht]
\centering
\begin{tabular}{||c c c c c c c c c||}
\hline
Case & $N_{\text{int}}$ & $N_{\text{sb}}$ & $K-1$ & $\bar{d}$ & $\lambda$ & $\mathrm{E}_{T}$ & $\|I - I^{\ast}\|_{L^{2}}$ & Training Time (sec.) \\ [0.5ex] 
\hline
1 & 8192 & 4096 & 4 & 20 & 0.1 & 0.0001 & 4.24e-05 & 23 \\ 
\hline
2 & 8192 & 4096 & 4 & 20 & 1 & 0.0002 & 3.77e-05 & 16 \\ 
\hline
3 & 8192 & 4096 & 4 & 20 & 0.1 & 0.0008 & 4.09e-05 & 30 \\ 
\hline
\end{tabular}
\caption{Results for the 1D infinite wall with a Gaussian source term.}
\label{table_2}
\end{table}
\subsubsection{Infinite wall with discontinuous source term}
In this section, we have presented a new model for the radiation transfer equation in graded-index media. Several types of deficiencies arise with traditional methods like the LSFEM and generalized  GFEM, as noted in \cite{zhao2012deficiency}. 
This model investigates heat transfer through radiation in an infinite slab with discontinuous source term. For this scenario, the exact analytical solution is provided as
\begin{equation}
I(x, \mu) = \mathbf{1}_{(x-0.5L)}\exp\left(-k_e(x-0.5L)\right) + \mathbf{1}_{(0.5L -x)}.
\end{equation}
In this equation, the constant parameter \(k_{e}\) indicates the strength of extinction, while \(L\) represents the thickness of the slab. The unit step function is denoted by \(\mathbf{1}_{[.]}\). The medium is homogeneous, and the walls are considered black. The source term can easily be derived from an exact solution. 
 
\begin{figure}[htbp]
\centering
\includegraphics[height=0.3\textheight]{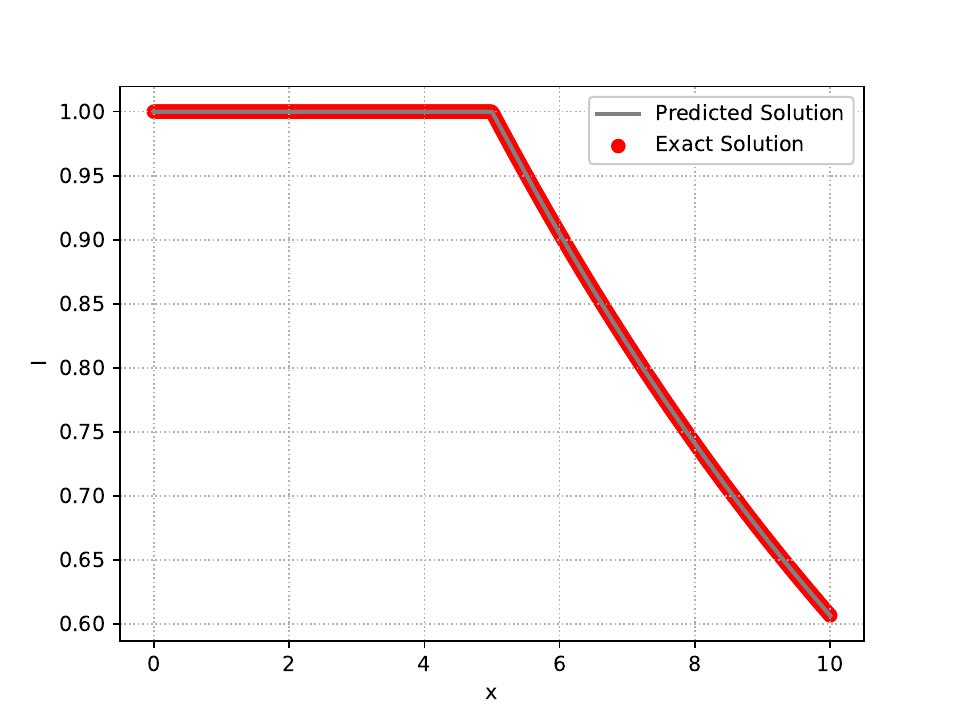}
\caption{Radiation intensity distributions for an infinite wall with a discontinuous source, comparing results from the exact solution and the PINN method at \( k_e = 0.1m^{-1}\).}
\label{2ka0.1}
\end{figure}

\begin{figure}[htbp]
\centering
\includegraphics[height=0.3\textheight]{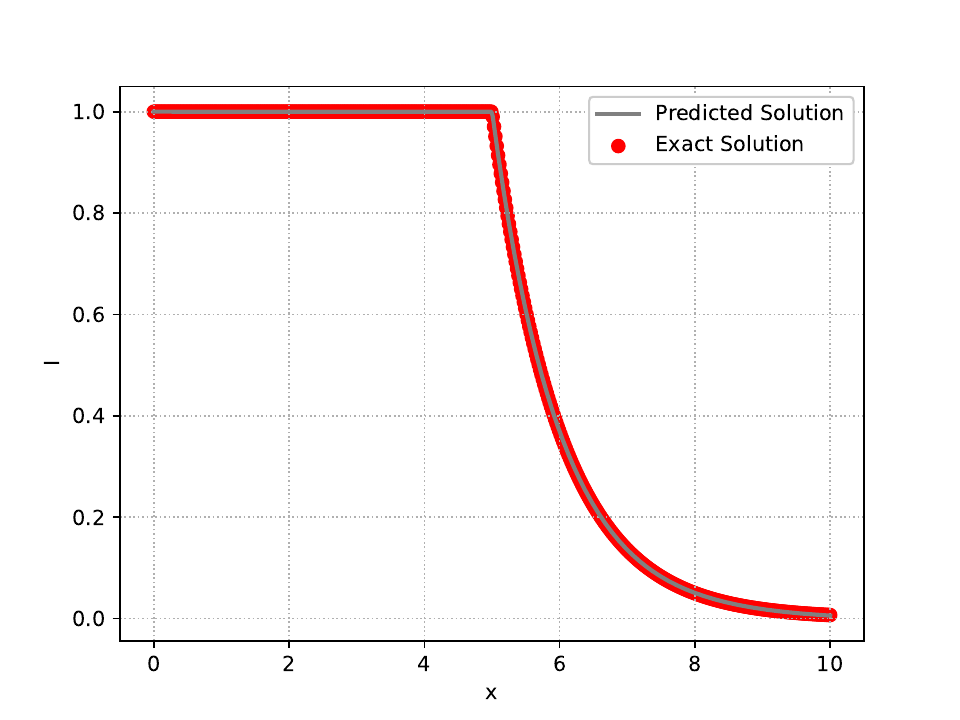}
\caption{Radiation intensity distributions for an infinite wall with a discontinuous source, comparing results from the exact solution and the PINN method at an extinction coefficient of \( k_e = 1m^{-1}\).}
\label{2ka1}
\end{figure}

\begin{figure}[htbp]
\centering
\includegraphics[height=0.3\textheight]{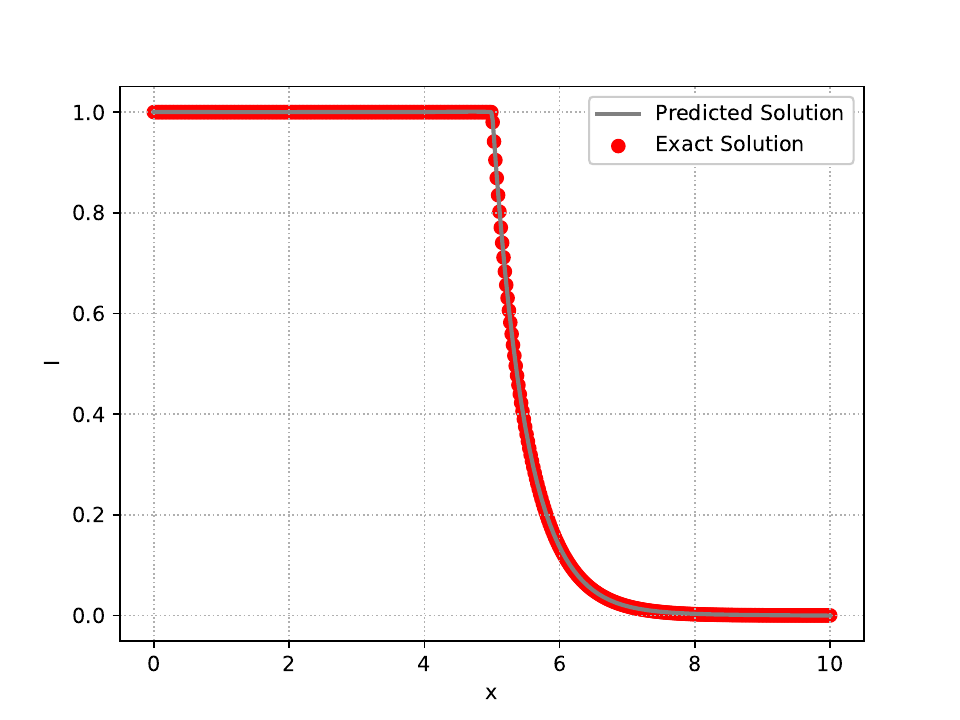}
\caption{Radiation intensity distributions for an infinite wall with a discontinuous source, comparing results from the exact solution and the PINN method at an extinction coefficient  of \( k_e = 2m^{-1}\).}
\label{2ka2}
\end{figure}

\begin{figure}[htbp]
\centering
\includegraphics[height=0.3\textheight]{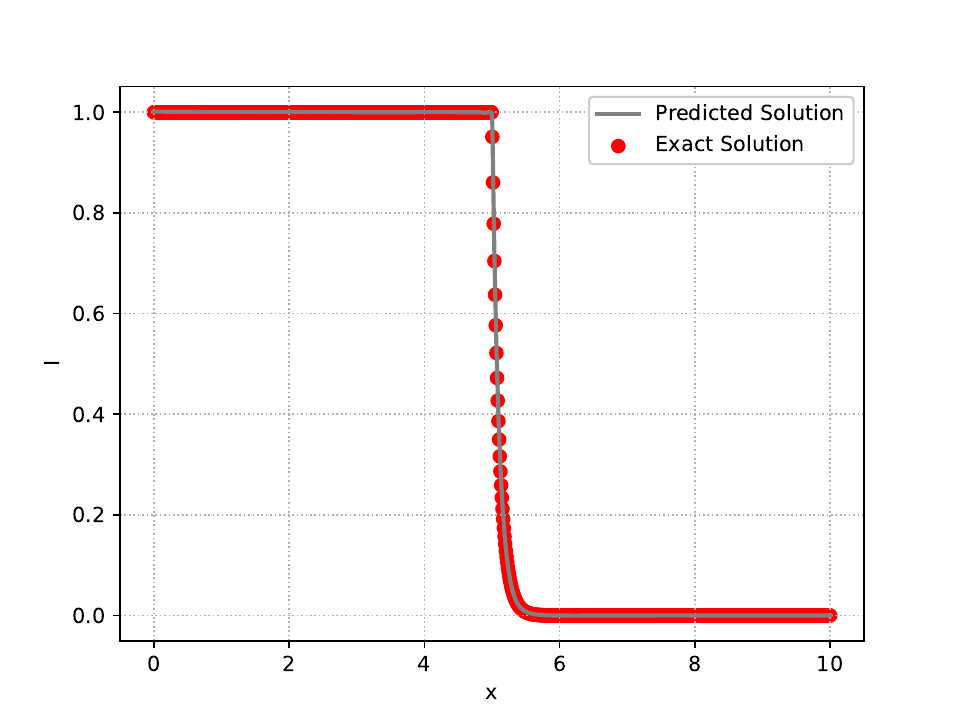}
\caption{Radiation intensity distributions for an infinite wall with a discontinuous source, comparing results from the exact solution and the PINN method at an extinction coefficient of \( k_e = 10m^{-1}\).}
\label{2ka10}
\end{figure}

The radiation intensity distributions solved for the direction \( \mu = 1 \) and \( L = 10 \) are compared across media with varying extinction coefficients: \( k_e = 0.1m^{-1} \), \( k_e = 1 \, m^{-1} \), \( k_e = 2m^{-1} \), and \( k_e = 10m^{-1} \). The boundary conditions of walls are unity.
The predicted and exact solutions of the RTE are illustrated in Figures \ref{2ka0.1}, \ref{2ka1}, \ref{2ka2}, and \ref{2ka10} for \( k_e = 0.1m^{-1} \), \( k_e = 1m^{-1} \), \( k_e = 2m^{-1} \), and \( k_e = 10m^{-1} \), respectively. The results obtained from the PINN demonstrate stability and agreement with the exact solutions across all four scenarios. Table \ref{table_3} shows that errors remain low at both boundaries, demonstrating the PINN's capability to approximate the PDE with minimal computational effort accurately.

\begin{table}[!ht]
\centering
\begin{tabular}{||c c c c c c c c c||}
\hline
Case & $N_{\text{int}}$ & $N_{\text{sb}}$ & $K-1$ & $\bar{d}$ & $\lambda$ & $\mathrm{E}_{T}$ & $\|I - I^{\ast}\|_{L^{2}}$ & Training Time (sec.) \\ [0.5ex] 
\hline
1 & 8192 & 4096 & 8 & 24 & 0.1 & 0.0005 & 0.00010 & 22 \\ 
\hline
2 & 8192 & 4096 & 8 & 24 & 0.1 & 0.0012 & 0.0004 & 18 \\ 
\hline
3 & 8192 & 4096 & 8 & 24 & 0.1 & 0.0020 & 0.0005 & 19 \\ 
\hline
4 & 8192 & 4096 & 8 & 24 & 1 & 0.009 & 0.0009 & 26 \\ 
\hline
\end{tabular}
\caption{Results for the radiation transport equation with an infinite wall with discontinuous source.}
\label{table_3}
\end{table}

\subsubsection{Square enclosure radiative distribution with discontinuous source along the diagonal}
\begin{figure}[htbp]
\centering
\includegraphics[height=0.4\textheight]{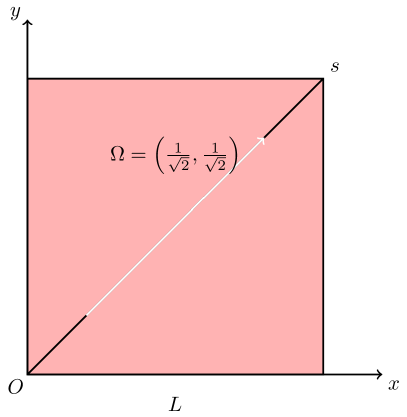}
\caption{Schematic of the square solution domain.}
\label{dir1}
\end{figure}

\begin{figure}[htbp]
\centering
\includegraphics[height=0.3\textheight]{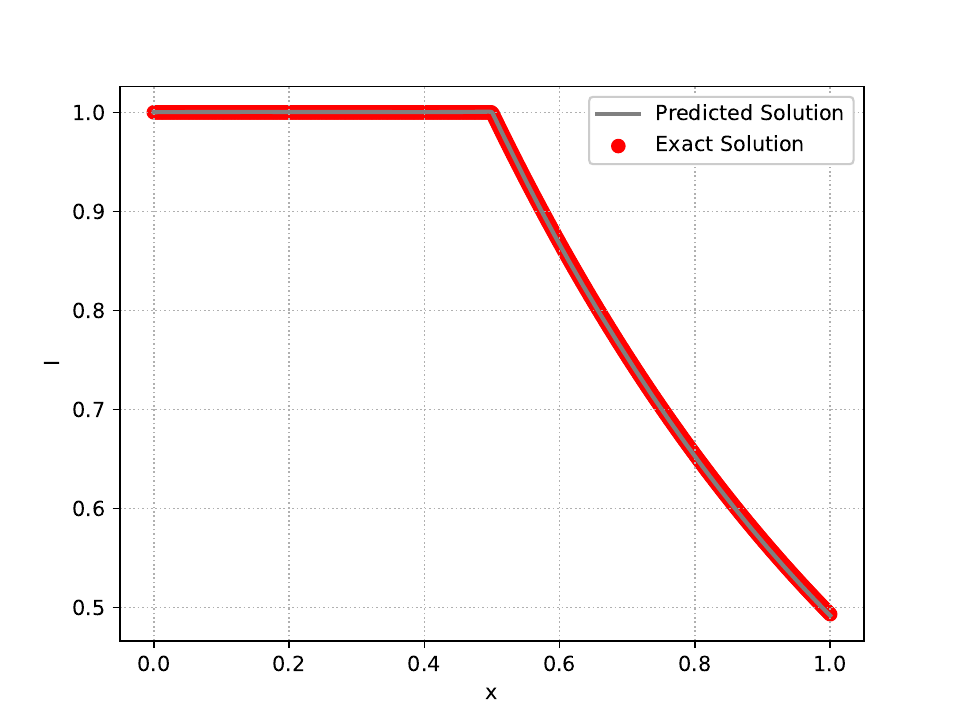}
\caption{Radiation distribution along the diagonal  at  extinction coefficients  \( k_e = 1m^{-1}\).}
\label{3ka1}
\end{figure}

\begin{figure}[htbp]
\centering
\includegraphics[height=0.3\textheight]{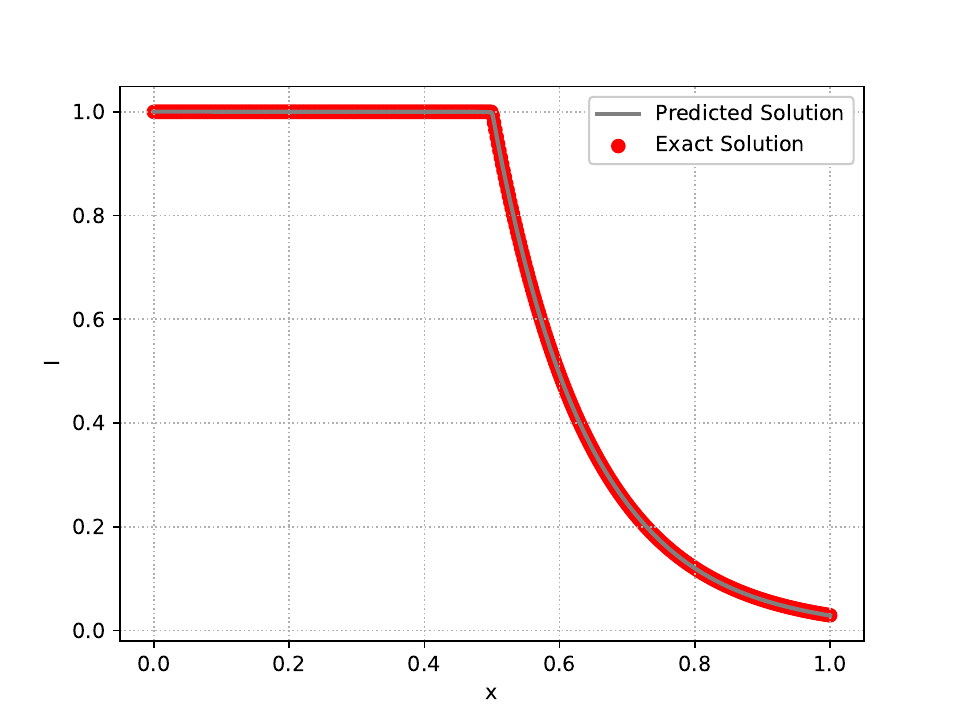}
\caption{Radiation distribution along the diagonal  at  extinction coefficients \( k_e = 5m^{-1}\).}
\label{3ka5}
\end{figure}

\begin{figure}[htbp]
\centering
\includegraphics[height=0.3\textheight]{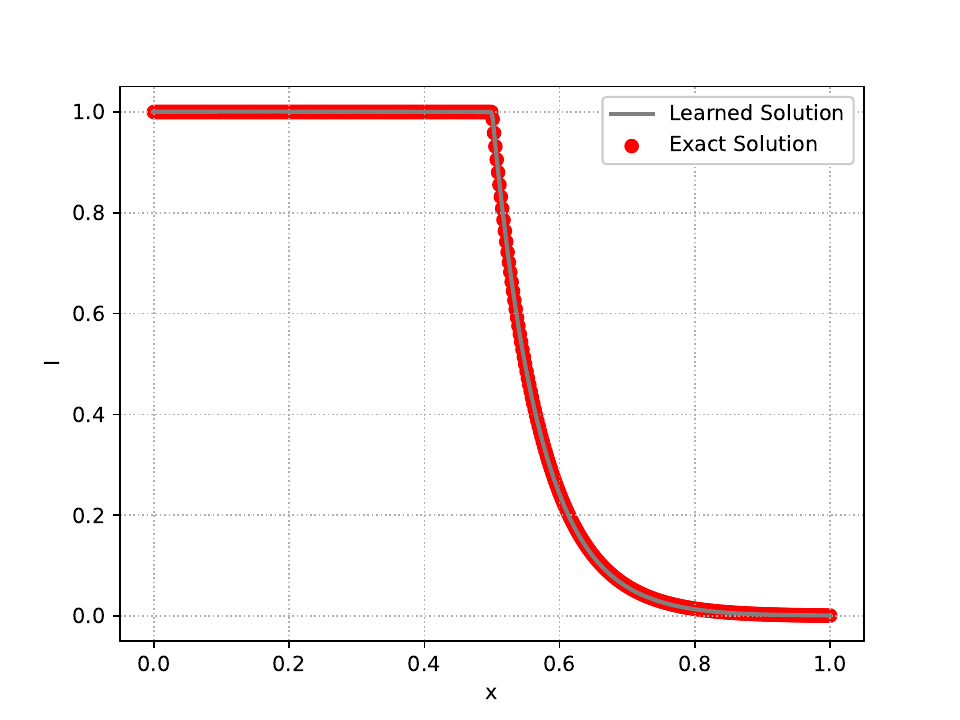}
\caption{Radiation distribution along the diagonal  at  extinction coefficients \( k_e = 10m^{-1}\).}
\label{3ka10}
\end{figure}

\begin{figure}[htbp]
\centering
\includegraphics[height=0.35\textheight]{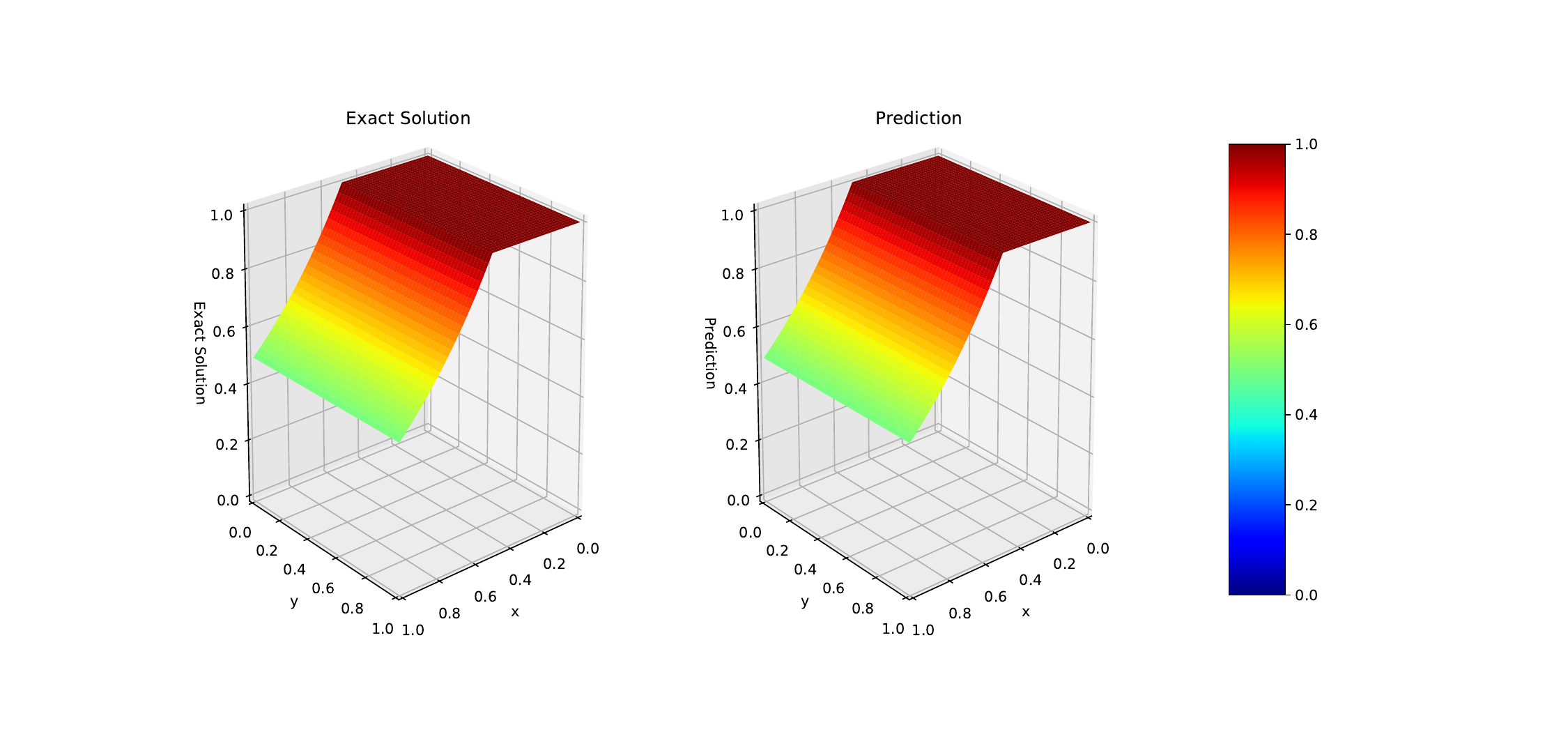}
\caption{3D plot for radiation distribution along the diagonal at extinction coefficients \( k_e = 1m^{-1}\).}
\label{3ka1a}
\end{figure}

\begin{figure}[htbp]
\centering
\includegraphics[height=0.35\textheight]{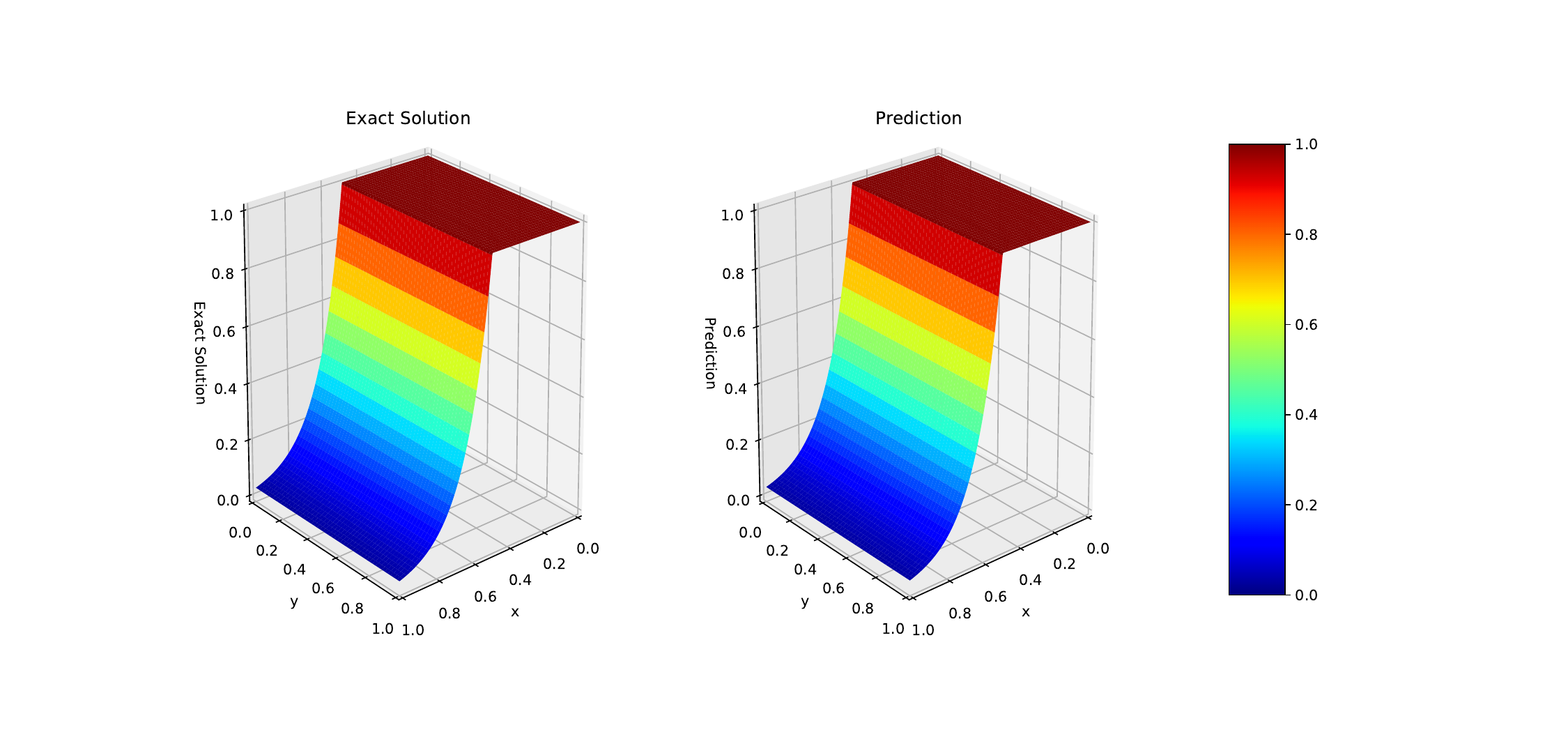}
\caption{3D plot for radiation distribution along the diagonal at  extinction coefficients \( k_e = 5m^{-1}\).}
\label{3ka5a}
\end{figure}

\begin{figure}[htbp]
\centering
\includegraphics[height=0.35\textheight]{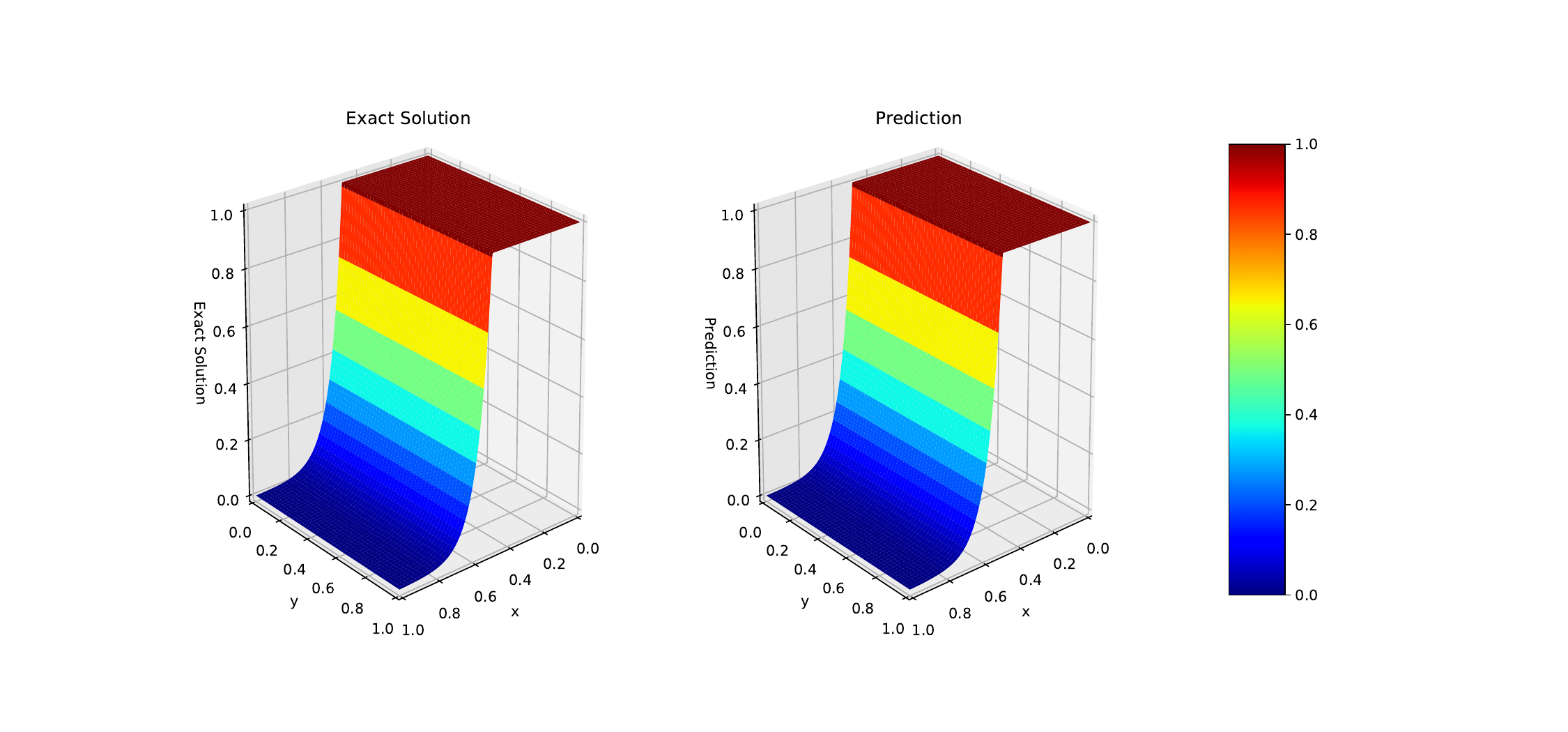}
\caption{3D plot for radiation distribution along the diagonal at  extinction coefficients \( k_e = 10m^{-1}\).}
\label{3ka10a}
\end{figure}

This section presents a new model for the radiation transfer equation in a square enclosure.
A 2D RTE problem is analyzed in this scenario. Figure \ref{dir1} shows the layout of the square medium. An exact solution exists for this configuration is
\begin{equation}
I(x,y,\mu) = \mathbf{1}_{(x+y-L)}\exp{\left[ -k_{e}\frac{x+y-L}{\sqrt{2} }\right] }+ \mathbf{1}_{(L-x-y)}, ~~~~ x, y \in [0,L]^2
\end{equation} 
In this equation, \(k_{e}\) denotes the extinction strength, \(1_{[.]}\) represents the unit step function. Let the direction \( \boldsymbol{\Omega} \) be \( \left[ \frac{1}{\sqrt{2}}, \frac{1}{\sqrt{2}} \right] \) for an enclosure with a side length of \( L = 1 \). The left and bottom walls have boundary conditions with an intensity of one. The source term can easily be derived from an exact solution.  Figures \ref{3ka1}, \ref{3ka5}, and \ref{3ka10} present the predicted and exact solutions of the RTE for \( k_e = 1 m^{-1} \), \( k_e = 5m^{-1} \), and \( k_e = 10 m^{-1} \), respectively. As illustrated in Table \ref{table_4}, the errors remain negligible at both boundaries, further showcasing the PINN's ability to accurately model the PDE while maintaining low computational expenses. We also present 3D plots in Figures \ref{3ka1a}, \ref{3ka5a}, and \ref{3ka10a}, which display the exact and predicted solutions of the RTE for \( k_e = 1m^{-1} \), \( k_e = 5m^{-1}\), and \( k_e = 10m^{-1}  \), respectively.
\begin{table}[!ht]
\centering
\begin{tabular}{||c c c c c c c c c||}
\hline
Case & $N_{\text{int}}$ & $N_{\text{sb}}$ & $K-1$ & $\bar{d}$ & $\lambda$ & $\mathrm{E}_{T}$ & $\|I - I^{\ast}\|_{L^{2}}$ & Training Time (sec.) \\ [0.5ex] 
\hline
1 & 8192 & 4096 & 8 & 28 & 1 & 0.007 & 0.00056 & 18 \\ 
\hline
2 & 8192 & 4096 & 8 & 28 & 1 & 0.0098 & 0.00048 & 19 \\ 
\hline
3 & 8192 & 4096 & 8 & 28 & 1 & 0.020 & 0.00059 & 20 \\ 
\hline
\end{tabular}
\caption{Results for square medium radiative transfer along the diagonal.}
\label{table_4}
\end{table}

\subsubsection{2D radiation distribution with Gaussian source term}\label{sec:square_gaussian_emissive_field}
\begin{figure}[htbp]
\centering
\includegraphics[height=0.35 \textheight]{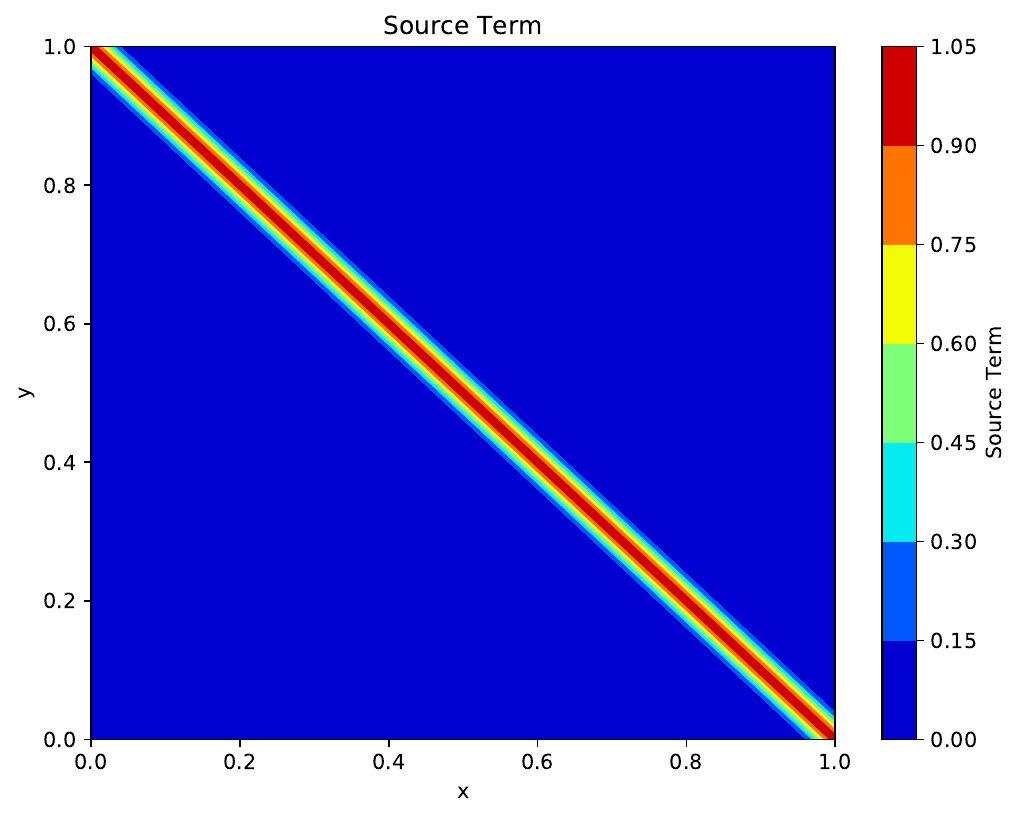}
\caption{Source term at $x$ is not equal to $y$.}
\label{source3}
\end{figure}

\begin{figure}[htbp]
\centering
\includegraphics[height=0.3\textheight]{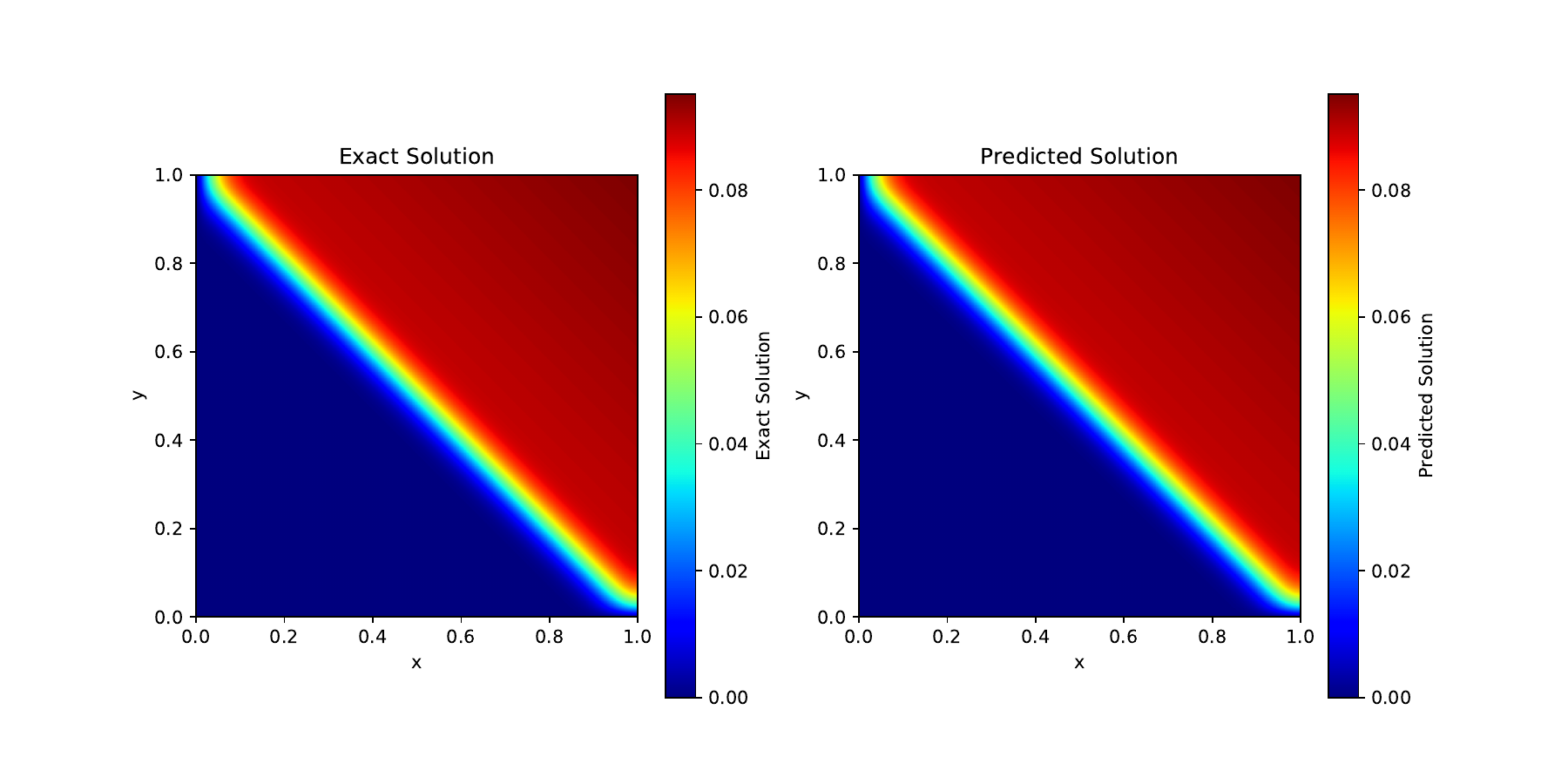}
\caption{2D Radiation distribution with a Gaussian source, solved by PINN and exact at \( k_e = 0.1m^{-1}\).}
\label{6ka0.1a}
\end{figure}

\begin{figure}[htbp]
\centering
\includegraphics[height=0.3\textheight]{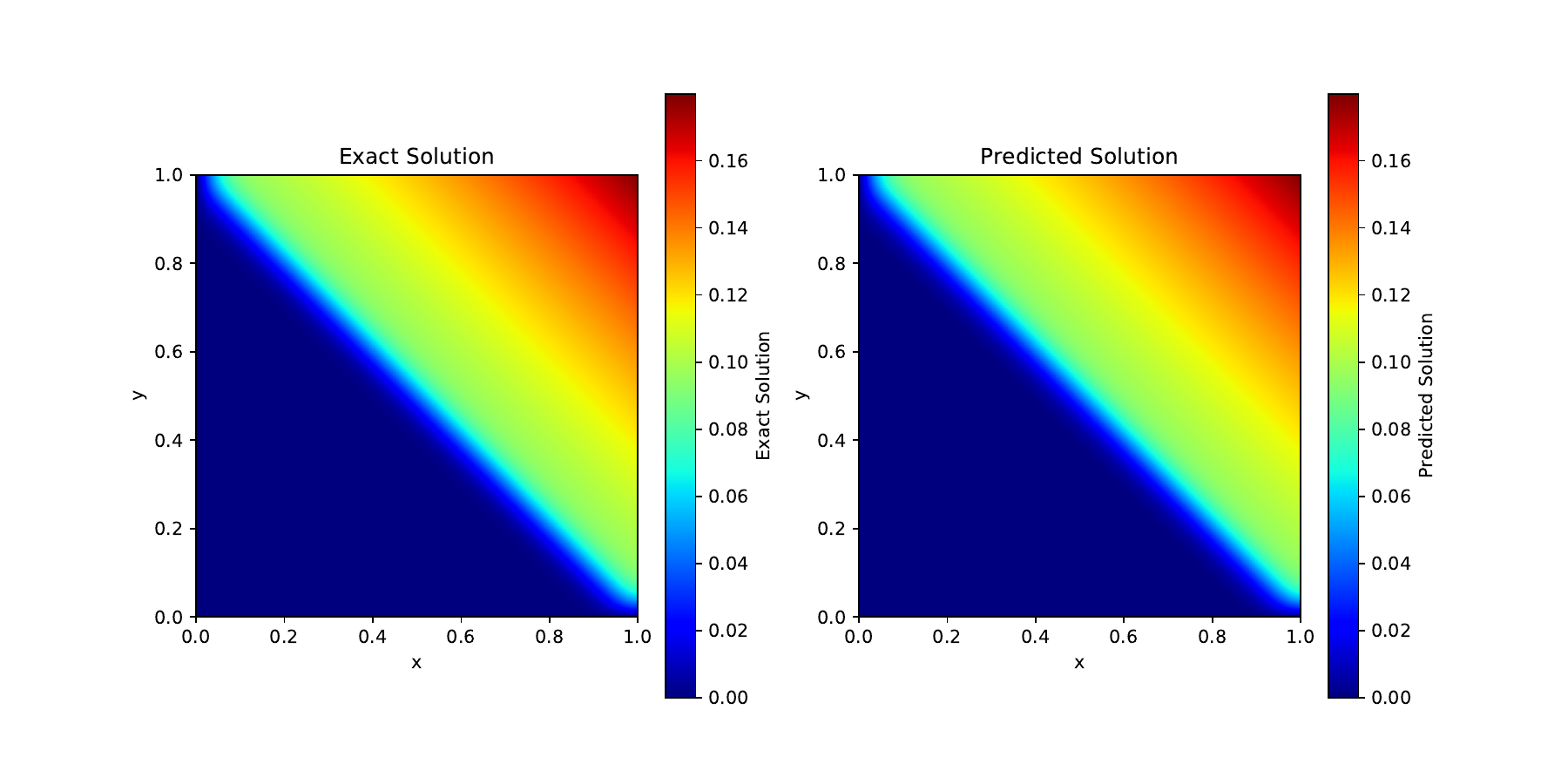}
\caption{2D Radiation distribution with a Gaussian source, solved by PINN and exact at  \( k_e = 1m^{-1}\).}
\label{6ka1}
\end{figure}

This model is based on an exponential square enclosure. This test represents the 2D form of Eq. (\ref{eq:1d}). The radiation transport equation can be modeled as
\begin{equation}\label{eq:2D}
\mu \frac{dI}{dx}+ \eta\frac{dI}{dy}+ k_{e}I = \exp\left(-\frac{\left( \frac{x+y}{\sqrt{2}}-c\right)^2}{\alpha^{2}}\right),  ~~~ x,y \in [0,1]
\end{equation}
Let the incident direction be defined as $\mu =\frac{\sqrt{2}}{2}$ , $\eta = \frac{\sqrt{2}}{2}$ and $ c= \frac{\sqrt{2}}{2}$.
The analytical solution \cite{zhao2013second} is
\begin{equation}
\begin{split}
I(x,y) &= \frac{\alpha \sqrt{\pi}}{2}\exp\left\lbrace -k_e\left[ \frac{|x+y|-|x-y|}{\sqrt{2}} -\frac{1-|x-y|}{\sqrt{2}} - \frac{\alpha^{2}k_{e}}{4}\right] \right\rbrace \\
&\times \left[ \operatorname{erf}\left( \frac{\alpha k_{e}}{2} + \frac{1-|x-y|}{\sqrt{2}\alpha}\right) -\operatorname{erf}\left( \frac{\alpha k_{e}}{2} + \dfrac{\frac{|x+y|-|x-y|}{\sqrt{2}} -\frac{1-|x-y|}{\sqrt{2}}} {\alpha}\right)  \right] 
\end{split}
\end{equation}
Fig. \ref{source3} illustrates the source term in Eq.(\ref{eq:2D}). Figs. \ref{6ka0.1a}, and \ref{6ka1} depict both the exact and predicted solutions of the RTE for \( k_e = 0.1m^{-1} \) and \( k_e = 1m^{-1} \) respectively. The errors, as presented in Table \ref{table_7}, are minimal, further highlighting the ability of the PINN with minimal computational effort accurately. At \(k_e = 1\,\mathrm{m}^{-1}\), the relative average error of the meshfree method \cite{zhao2013second} is \(0.7\). The relative error of the MRT lattice Boltzmann method \cite{feng2021performance} at \( k_e = 1m^{-1}\) is $2.5\%$. \cite{yi2016lattice}  simulated the same model using the lattice Boltzmann method and compared the results with the GFEM. Both \cite{yi2016lattice} and \cite{zhao2013second} showed that for small values of \(k_{e}\), the GFEM and the meshfree method, respectively, performed poorly.

\begin{table}[!ht]
\centering
\begin{tabular}{||c c c c c c c c c||}
\hline
Case & $N_{\text{int}}$ & $N_{\text{sb}}$ & $K-1$ & $\bar{d}$ & $\lambda$ & $\mathrm{E}_{T}$ & $\|I - I^{\ast}\|_{L^{2}}$ & Training Time (sec) \\ [0.5ex] 
\hline
1 & 8192 & 4096 & 4 & 20 & 0.1 & 0.0008 & 0.04 & 31 \\ 
\hline
2 & 8192 & 4096 & 4 & 20 & 0.1 & 0.0009 & 0.05 & 27 \\ 
\hline
\end{tabular}
\caption{Results for 2D Radiation distribution with a Gaussian source(forword).}
\label{table_7}
\end{table}
\newpage
\subsubsection{Radiation distribution with a Gaussian source field along diagonal}
\begin{figure}[htbp] \centering \includegraphics[height=0.35\textheight]{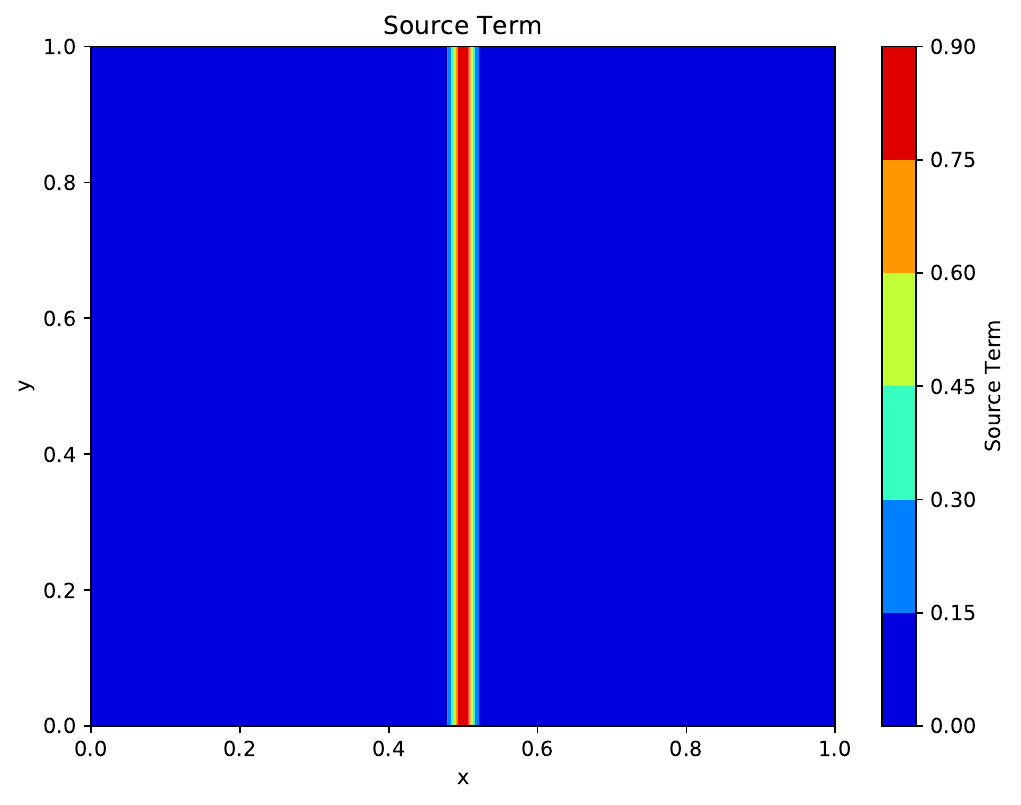} \caption{Source term at $x = y$.} \label{source2} \end{figure}

\begin{figure}[htbp] \centering \includegraphics[height=0.3\textheight]{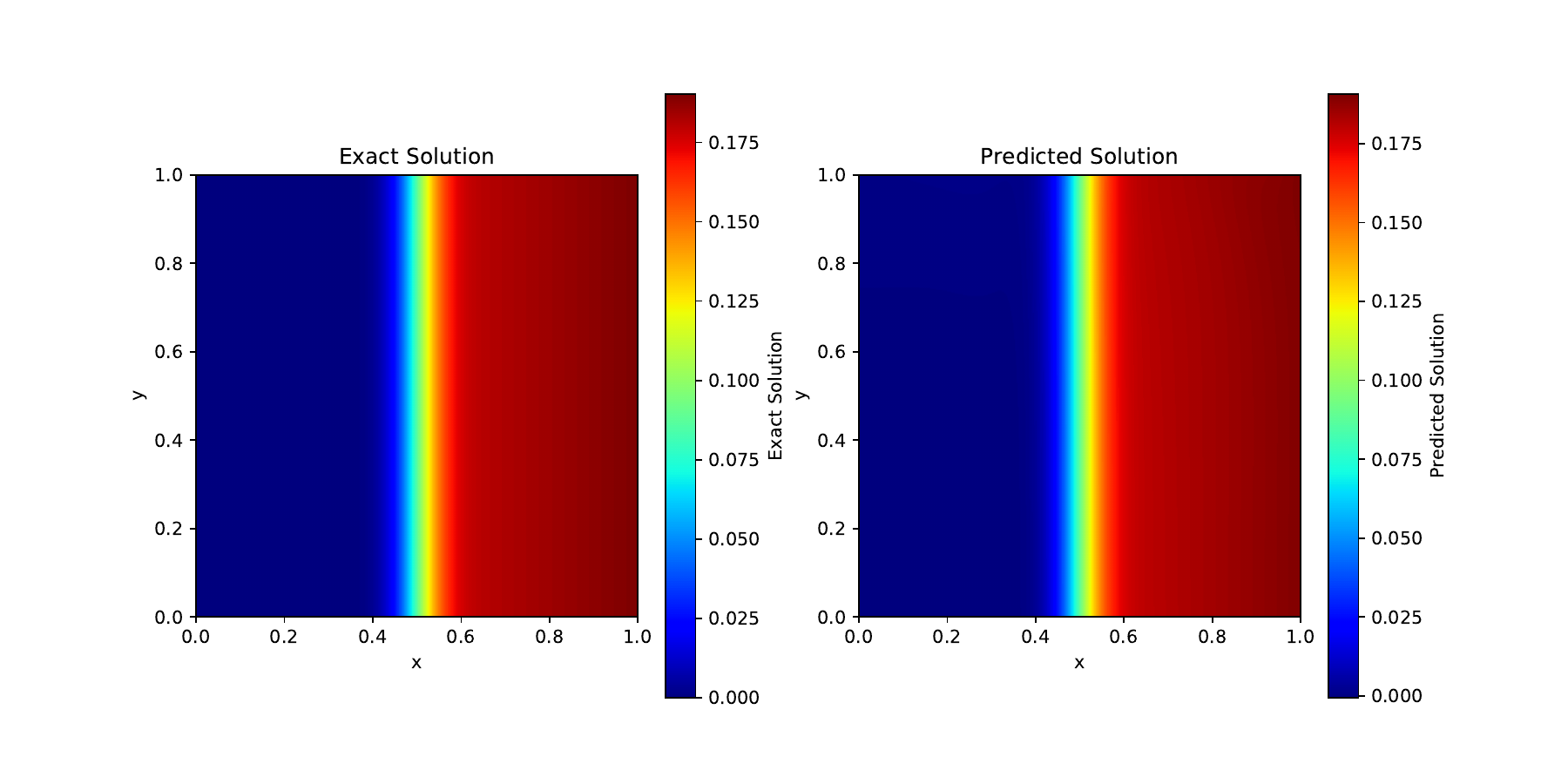} \caption{Radiation distribution with a Gaussian source, solved by PINN and exact at \( k_e = 0.1m^{-1}\) for $x = y$.} \label{5ka.1} \end{figure}

\begin{figure}[htbp] \centering \includegraphics[height=0.3\textheight]{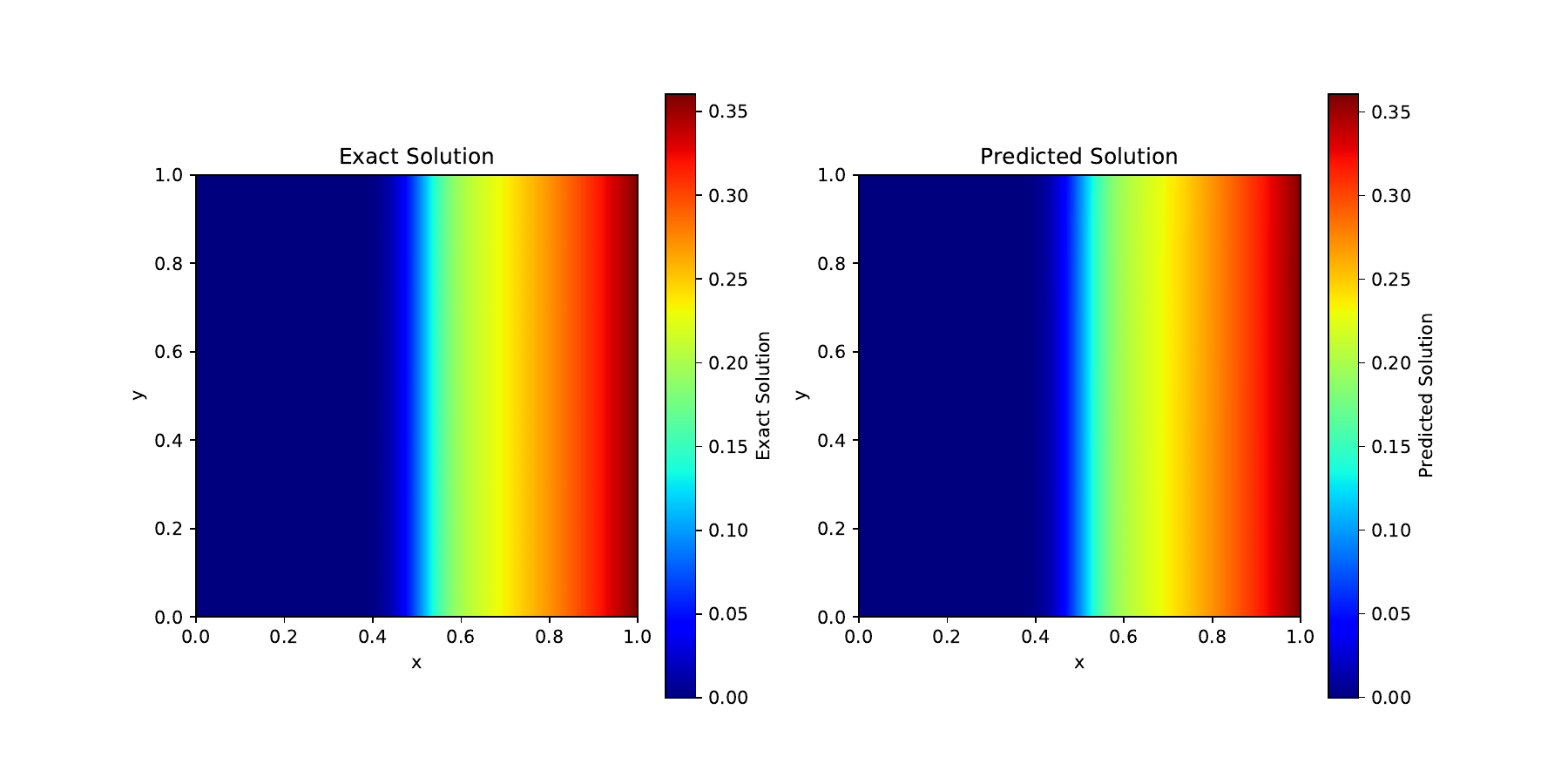} 
\caption{Radiation distribution with a Gaussian source, solved by PINN and exact at \( k_e = 1m^{-1}\) for $x = y$.} \label{5ka1} \end{figure}

\begin{figure}[htbp] \centering \includegraphics[height=0.3\textheight]{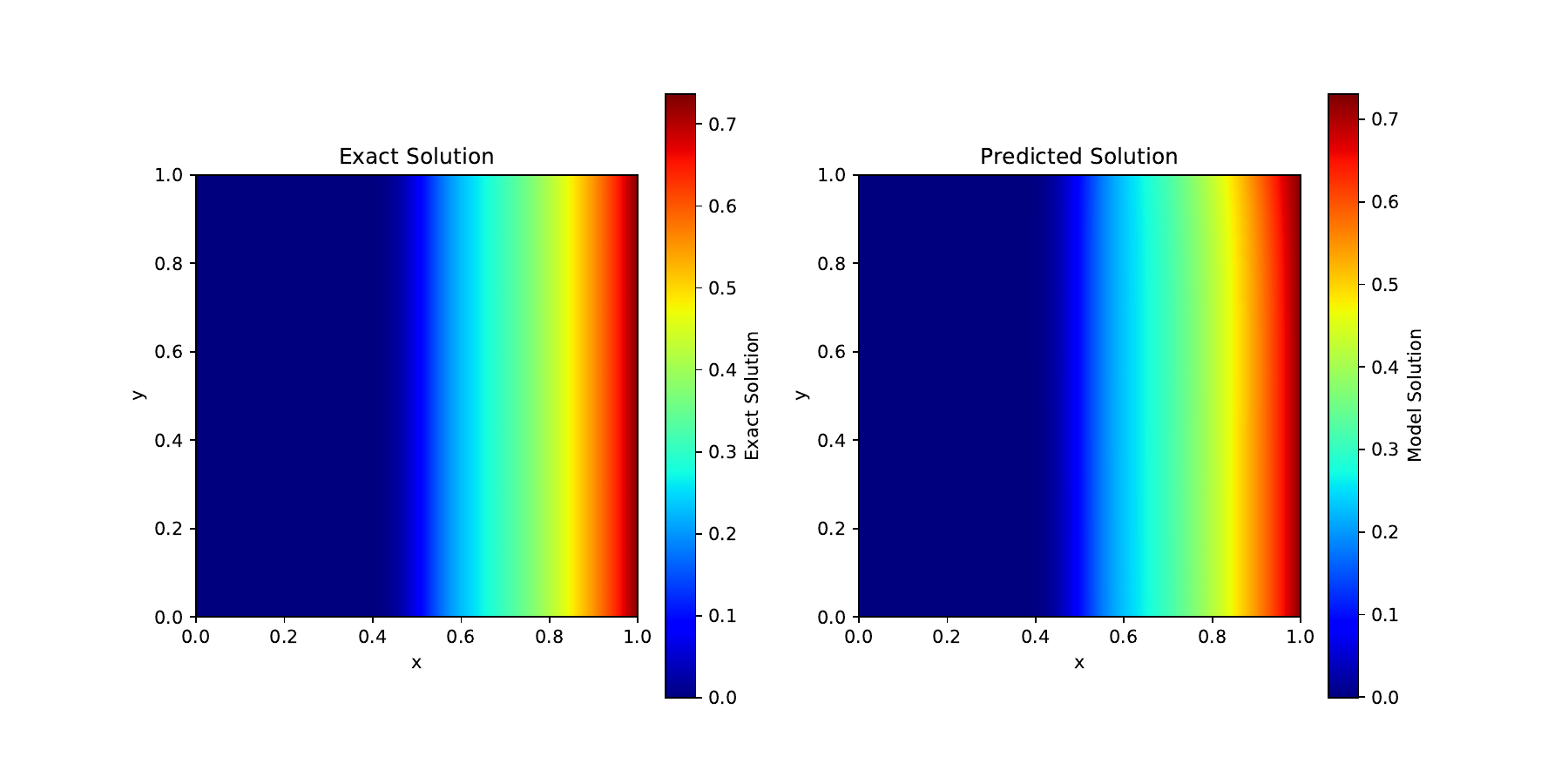} \caption{Radiation distribution with a Gaussian source, solved by PINN and exact at \( k_e = 2m^{-1}\) for $x = y$.} \label{5ka2} \end{figure}

The numerical experiment discussed in Section~(\ref{sec:square_gaussian_emissive_field}) is conducted under the condition \(x = y\) with black and cold boundary conditions. A PINN determines the radiation distribution along the square enclosure's diagonal (\(y = x\)). Figure \ref{source2} represents the source term. Figures \ref{5ka.1}, \ref{5ka1}, and \ref{5ka2} illustrate exact and predicted solutions of the RTE for \( k_e = 0.1m^{-1} \), \( k_e = 1m^{-1} \), and \( k_e = 2m^{-1}\), respectively. As shown in Table \ref{table_6}, the errors are minimal at both boundaries, further demonstrating the PINN's capability to approximate the PDE with low computational cost accurately.
% Surface plots of the exact and learned solutions for $k_a = 0.1$, $k_a = 1$, and $k_a = 2$ are shown in figures \ref{4ka0.1a}, \ref{4ka1a}, and \ref{4ka2a}.

\begin{table}[!ht]
\centering
\begin{tabular}{||c c c c c c c c c||}
\hline
Case & $N_{\text{int}}$ & $N_{\text{sb}}$ & $K-1$ & $\bar{d}$ & $\lambda$ & $\mathrm{E}_{T}$ & $\|I - I^{\ast}\|_{L^{2}}$ & Training Time (sec.) \\ [0.5ex] 
\hline
1 & 8192 & 4096 & 4 & 20 & 0.1 & 0.00029  & 0.0002 & 9 \\ 
\hline
2 & 8192 & 4096 & 4 & 20 & 0.1 & 0.0008 & 0.0012 & 10 \\ 
\hline
3 & 8192 & 4096 & 4 & 20 & 0.1 & 0.0002 & 0.003 & 17 \\ 
\hline
\end{tabular}
\caption{Results of the radiation transport equation for Gaussian-shaped emissive field along the diagonal.}
\label{table_6}
\end{table}
\subsection{Inverse problems}
\subsubsection{2D radiation distribution with Gaussian source term}
\begin{figure}[htbp]
\centering
\includegraphics[height=0.3\textheight]{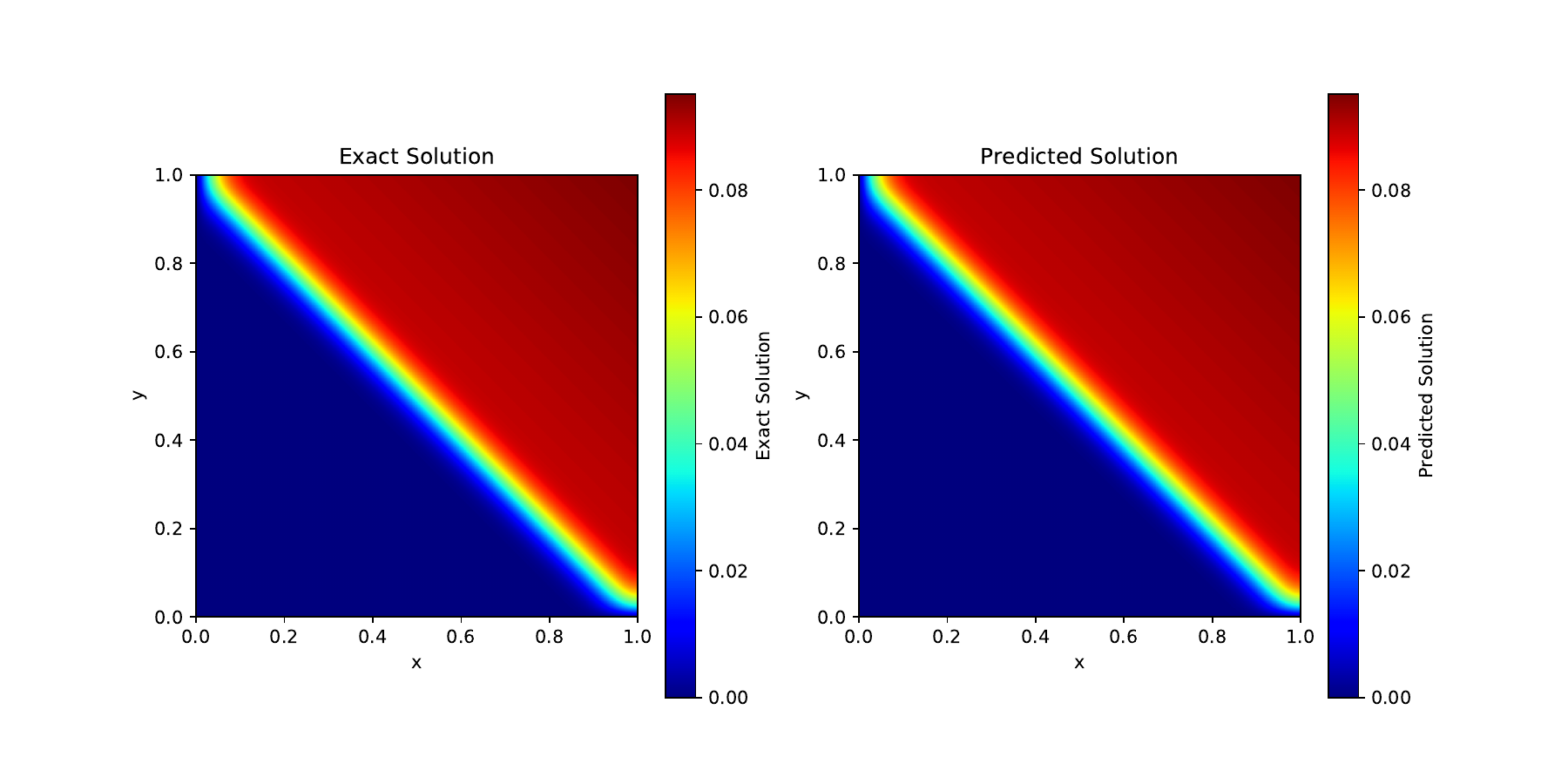}
\caption{2D Radiation distribution with a Gaussian source, solved by inverse PINN and exact at \( k_e = 0.1m^{-1}\).}
\label{7ka2a}
\end{figure}

\begin{figure}[htbp]
\centering
\includegraphics[height=0.3\textheight]{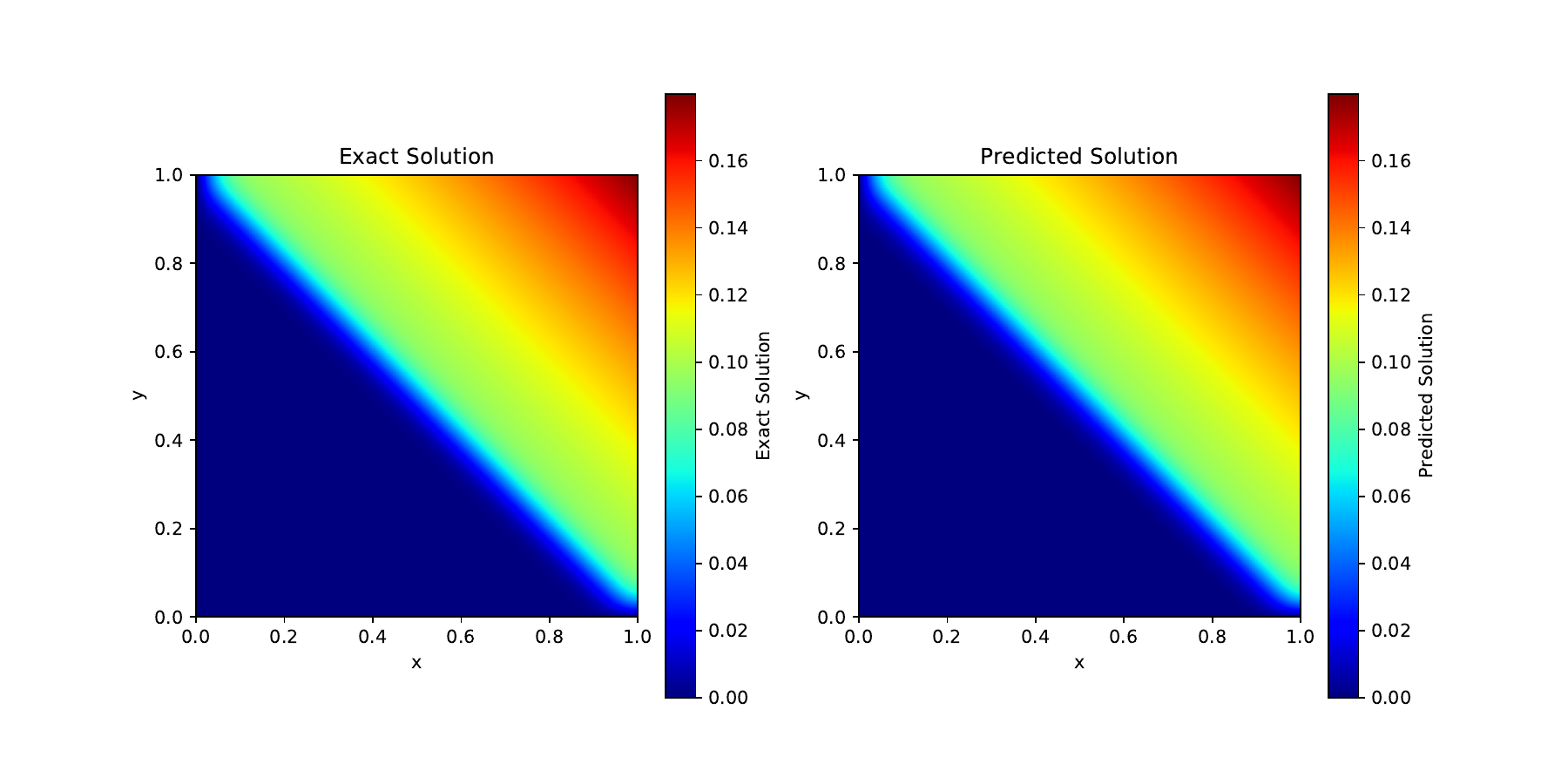}
\caption{2D Radiation distribution with a Gaussian source, solved by inverse PINN and exact at \( k_e = 1m^{-1}\).}
\label{7ka2aa}
\end{figure}

The numerical experiment of section (\ref{sec:square_gaussian_emissive_field}) is performed as an inverse problem using Algorithm \ref{alg2}. In this experiment, we excluded the boundary conditions. Figs.\ref{7ka2a} and \ref{7ka2aa} show the exact and predicted solutions of the RTE for \(k_e = 0.1m^{-1}\) and \(k_e = 1m^{-1}\), respectively. The errors, shown in Table \ref{table_8}, are minimal for different $k_{e}$ values, demonstrating that the PINN is with high accuracy and at a meager computational cost. At \(k_e = 1\,\mathrm{m}^{-1}\), the relative average error of the mass-free method \cite{zhao2013second} is \(0.7\). At \(k_e = 1\,\mathrm{m}^{-1}\), the relative average error of the mass-free method(MSORTE form) \cite{zhao2013second} is \(0.04\).
 At $k_{e} = 1m^{-1}$, the relative error of the MRT lattice Boltzmann algorithm, as simulated by Feng et al. \cite{feng2021performance}, is $2.5\%$. The PINN method demonstrates a relative error of \(0.3\%\) at \(k_{e} = 1\), notably lower than that of other methods.

\begin{table}[!ht]
\centering
\begin{tabular}{||c c c c c c c c c||}
\hline
Case & $N_{\text{int}}$ & $N_{\text{d}}$ & $K-1$ & $\bar{d}$ & $\lambda$ & $\mathrm{E}_{T}$ & $\|I - I^{\ast}\|_{L^{2}}$ & Training Time (sec.) \\ [0.5ex] 
\hline
1 & 16384 & 8192 & 4 & 20 & 0.1 & 0.003 & 0.0008 & 36 \\ 
\hline
2 & 16384 & 8192 & 4 & 20 & 0.1 & 0.003 & 0.0005 & 76 \\ 
\hline

\end{tabular}
\caption{Results for the 2D Radiation distribution with a Gaussian source(inverse).}
\label{table_8}
\end{table}
\section{Conclusion} \label{sec:4}
Solving the RTE with graded index (\ref{eqn:A}) presents a formidable challenge due to their inherent high-dimensional nature, especially when considering the most general scenario with high dimensions. Furthermore, incorporating diverse physical phenomena such as emission, absorption, and scattering, alongside the variability in optical parameters across the medium, adds complexity to devising efficient numerical algorithms. As highlighted earlier, existing methodologies often grapple with the curse of dimensionality, necessitating substantial computational resources to attain the desired level of precision. In our study, we propose an innovative solution to this challenge. Our approach, outlined in Algorithms \ref{alg1} and \ref{alg2}, harnesses PINNs and sophisticated neural network architectures tailored specifically for approximating the radiative intensity outlined in the equation. Through iterative training using gradient descent, our network endeavors to minimize comprehensive loss function (\ref{eqn:La}) and (\ref{eqn:Lb}), respectively. This loss function encapsulates the residual error from integrating the neural network representation into the RTE (\ref{eqn:A}). The residual errors are strategically positioned at training points, aligning with quadrature points based on a predefined quadrature rule. To alleviate the computational burden associated with high dimensionality, we employ Sobol low-discrepancy sequences as training points, optimizing efficiency while preserving accuracy.\\

Traditional numerical methods often encounter limitations when solving first-order equations, particularly in singular boundary conditions. To address these challenges, these methods typically convert the equations into second-order form to eliminate singularities at the boundaries before initiating the simulation. However, this reformulation can introduce inaccuracies, especially in scenarios with steep gradients or discontinuities. In their original first-order form, traditional methods are prone to numerical instabilities, often resulting in oscillations or wiggles in the solution. PINNs provide a robust alternative to overcome these challenges. PINNs embed the governing physics into the loss function by directly solving the first-order equations without requiring reformulation. This approach allows them to naturally handle boundary conditions and produce smooth, stable solutions, even in cases where traditional methods fail.
Therefore, we posit that PINN algorithms \ref{alg1} and \ref{alg2} serve as a versatile, straightforward-to-implement, swift, and precise simulator for RTE phenomena. We demonstrate that this algorithm excels in speed and accuracy through numerical experiments. In essence, we contribute novel machine learning methodologies that offer a swift, user-friendly, and precise means of simulating various facets of radiation heat transfer in graded index phenomena. Our results showed that this approach worked well to reduce numerical errors, as there were no strange bumps in the results. 
\section*{Declaration of competing interest}The authors declare that they have no competing interests.
\section*{Acknowledgment}
The first author acknowledges the Ministry of Human Resource Development (MHRD), Government of India, for providing institutional funding and support at IIT Madras.
\section*{Appendix}
An estimate of the generalization error for equation (\ref{eqn:A}) is provided for the forward in \ref{thm:1}. An estimate of the generalization error for steady state equation (\ref{eqn:BB}) is provided for forward problems in \ref{thm:generalization_error_forward_inverse}.
\begin{appendixthm}\label{thm:1}
Assume \( I \in L^{2}(\boldsymbol{D}) \) is the unique weak solution to the RTE in graded-index media, where the coefficients \( 0 \leq k_{e}, k_{s} \in L^{\infty}(\textit{D}) \), and phase function \( \Phi(\boldsymbol{\Omega}, \boldsymbol{\Omega}') \in C^p(S\times S) \) for some \( p > 0 \). Let \( I^{*} = I_{\Theta} \in C^p(\boldsymbol{D}) \) be the solution generated by Algorithm \ref{alg1} (forward PINN applied to the equation). Assume the condition \( \max \left\lbrace HK_{V}(I^{\ast}),  HK_{V}(\mathrm{R}_{int, \Theta^{\ast}})\right\rbrace < \infty \), where \(  HK_{V} \) represents the Hardy-Krause variation.
Further, assume that Sobol points are used as training points \( \zeta_{int}, \zeta_{sb}, \zeta_{tb} \), and a Gauss-quadrature rule of order \( a = a(p) \) is applied to approximate the scattering term in residual. Additionally, we assume that a refractive index  \( n= f(s) \in \textit{D} \) are arbitrarily chosen in a manner consistent with the general model structure of the radiative transport equation.
Under these assumptions:
\[
n^2 k_a I_b(T_g) \leq n^2 k_a MI + f(s),
\]
\[
( \boldsymbol{\Omega} \cdot \nabla)I+ \frac{1}{n \sin \theta} \frac{\partial}{\partial \phi} \left( (s_1 \cdot \nabla n) I \right) + \frac{1}{n \sin \theta} \frac{\partial}{\partial \theta} \left( I(\boldsymbol{\Omega} \cos \theta - k) \cdot \nabla n \right) \leq nMI + ( \boldsymbol{\Omega} \cdot \nabla_{s})I, \quad M > 0.
\]
The estimation of the generalized error for the forward problem is

\begin{dmath}
(\mathrm{E}_{G})^{2} \leq V \left( (\mathrm{E}_{T}^{tb})^{2} + \nu (\mathrm{E}_{T}^{sb})^{2} + \nu (\mathrm{E}_{T}^{int})^{2} \right) + VV_{2} \left( \dfrac{(\log(N_{tb}))^{2d}}{N_{tb}} + \nu \dfrac{(\log(N_{sb}))^{2d}}{N_{sb}} + \nu \dfrac{(\log(N_{int}))^{2d+1}}{N_{int}} + \nu N_{\boldsymbol{S}}^{-2a} \right),
\end{dmath}
where:
\[
\bar{V} = (\|\boldsymbol{D}\|, \|\Phi\|_{C^p}, \|\hat{I}\|),
\]
\[
V = \left( T + \nu V_1 \bar{T}^2 \exp(\nu V_1 T) \right),
\]
\[
V_1 = \frac{2 \nu (\|k_s\|_{L^{\infty}} + \|\Sigma_g\|_{L^{\infty}})}{4\pi},
\]
\[
V_2 = \max \left\{  HK_{V}(\mathrm{R}^{\ast}_{tb})^2,  HK_{V}(\mathrm{R}^{\ast}_{sb})^2,  HK_{V}(\mathrm{R}^{\ast}_{int})^2, \bar{V} \right\}.
\]
\end{appendixthm}
\begin{proof}
We are following \cite{mishra2021physics}
\begin{dmath}
\mathscr{E}(I^{\ast}, \Phi) = \sum\limits_{i=1}^{N_{\boldsymbol{S}}}w_{i}^{\boldsymbol{S}}\Phi(\boldsymbol{\Omega},\boldsymbol{\Omega}_{i}^{\boldsymbol{S}}) I^{\ast}(t,s,\boldsymbol{\Omega}_{i}^{\boldsymbol{S}}) -\int\limits_{4\pi}\Phi(\boldsymbol{\Omega},\boldsymbol{\Omega}') I^{\ast}(t,s,\boldsymbol{\Omega}) d\Omega',
\end{dmath}
$\bar{I}=I^{\ast}-I,$
\begin{dmath}
\frac{n}{c_{0}}\frac{\partial}{\partial t}\bar{I} + (k_{e} + \boldsymbol{\Omega} \cdot \nabla)\bar{I} + \frac{1}{n\sin \theta}\frac{\partial}{\partial \theta}\left\lbrace \bar{I}(\boldsymbol{\Omega} \cos \theta - k) \cdot \nabla n \right\rbrace 
 +  \frac{1}{n\sin \theta}\frac{\partial}{\partial \phi}\left\lbrace (s_{1} \cdot \nabla n)\bar{I}\right\rbrace= n^{2}k_{a}\bar{I}_{b}(T_{g}) +\frac{k_{s}}{4\pi}\int\limits_{4\pi}\Phi(\boldsymbol{\Omega},\boldsymbol{\Omega}') \bar{I}(t,s,\boldsymbol{\Omega}) d\Omega'+ 
\mathscr{E}(I^{\ast}), \label{eq:eq1} \\
\end{dmath}
We define
\begin{align*}
\bar{I}(0, s,\boldsymbol{\Omega})&=\mathrm{R}^{\ast}_{tb},~~(s,\boldsymbol{\Omega}) \in \textit{D}\times S\\
\bar{I}(t, s,\boldsymbol{\Omega})&=\mathrm{R}^{\ast}_{sb}, ~~(t, s,\boldsymbol{\Omega}) \in \beta \label{eq:eq3}\tag{A}\\
\end{align*}
Multiplying  Eq.\eqref{eq:eq1} by $\bar{I}$ on both side,
\begin{dmath}
\frac{n}{2c_{0}}\frac{d \bar{I}^{2}}{dt}=- k_{e} \bar{I^{2}}-( \boldsymbol{\Omega} \cdot \nabla_{s})\left( \frac{\bar{I}^{2}}{2}\right) - \frac{1}{n\sin \theta}\frac{\partial}{\partial \theta}\left\lbrace \bar{I}(\boldsymbol{\Omega} \cos \theta - k) \cdot \nabla n \right\rbrace \bar{I}
 -  \frac{1}{n\sin \theta}\frac{\partial}{\partial \phi}\left\lbrace (s_{1} \cdot \nabla n)\bar{I}\right\rbrace\bar{I} + n^{2}k_{a}\bar{I}_{b}\bar{I}(T_{g})-\frac{k_{s}}{4\pi}\int\limits_{4\pi}\Phi(\boldsymbol{\Omega},\boldsymbol{\Omega}') \bar{I}(t, s,\boldsymbol{\Omega})\bar{I}((t, s,\boldsymbol{\Omega'})d{\Omega'}+ 
\mathscr{E}(I^{\ast}, \Phi)\bar{I}. \label{eq:eq4a}
\end{dmath}
We can observe that  $ ( \boldsymbol{\Omega} \cdot \nabla)I + \frac{1}{n\sin \theta}\frac{\partial}{\partial \phi}\left\lbrace (s_{1} \cdot \nabla n)\bar{I}\right\rbrace  + \frac{1}{n\sin \theta}\frac{\partial}{\partial \theta}\left\lbrace \bar{I}(\boldsymbol{\Omega} \cos \theta - k) \cdot \nabla n \right\rbrace  \leq M\bar{I} + ( \boldsymbol{\Omega} \cdot \nabla_{s})I$,\\
$  n^{2}k_{a}\bar{I}_{b}(T_{g}) \leq n^{2}k_{a}nM\bar{I} + f(s)  $,
\\
And
\begin{align}
\int\limits_{\textit{D}\times S} \vert f(s) \bar{I} \vert ds  d\boldsymbol{\Omega} \leq \int\limits_{\textit{D}\times S} \vert f\vert^{2}    ds  d\boldsymbol{\Omega}
+ \int\limits_{\textit{D}\times S} \vert \bar{I}\vert ^{2} ds  \boldsymbol{\Omega}.\label{eq:leto}
\end{align} 
We substitute the values into Eq.(\ref{eq:eq4a}), integrate the result over $\textit{D} \times S$, and apply integration by part
and Cauchy sequence and $k_{e}, k_{s}>0$, and $t \in (0,T]$
\begin{dmath}
\begin{split}
\frac{n}{2c_{0}}\frac{d}{dt}\int\limits_{\textit{D}\times S}\bar{I}^{2}(t, s,\boldsymbol{\Omega})ds  d\boldsymbol{\Omega}  &\leq 
\int\limits_{\textit{D}\times S}\bar{I}^{2}(t, s,\boldsymbol{\Omega})ds  d\boldsymbol{\Omega}
- \int\limits_{(\partial{\textit{D}} \times S)_{-}}(\Omega \cdot k(s)) \frac{\bar{I}^{2}(t, s,\boldsymbol{\Omega})}{2}dk(s) d\boldsymbol{\Omega}\\
& \quad+ \int\limits_{\textit{D}\times S} \frac{k_{s}}{4\pi}\int\limits_{4\pi}\Phi(\boldsymbol{\Omega},\boldsymbol{\Omega}') \bar{I}(t, s,\boldsymbol{\Omega})\bar{I}(t, s,\boldsymbol{\Omega}')d\Omega' d\boldsymbol{\Omega}  ds \\
& \quad + \int\limits_{\textit{D}\times S} \vert f\vert^{2}   ds  d\boldsymbol{\Omega} 
+ \int\limits_{\textit{D}\times S}\frac{(\mathscr{E}(I^{\ast},\Phi)(t, s,\boldsymbol{\Omega}))^{2}}{2}d\boldsymbol{\Omega}  ds \label{eq:letta}.
\end{split}
\end{dmath}
In this case, \( dk(s) \) signifies the surface measure on \( \partial \textit{D} \), and it can defined as
\[ \beta =(\partial{\textit{D}} \times S)_{-}=\left\lbrace (s,\boldsymbol{\Omega}) \in \partial{\textit{D}} \times S :\boldsymbol{\Omega}.k(s)\leq 0 \right\rbrace, \]  
And with  $k(s)$ unit normal at  $s \in \partial{\textit{D}}.$
Select a $\hat{T} \in (0,T]$ and integrate Eq.(\ref{eq:letta}) over  $ (0,\bar{T}),$
\begin{dmath}
\begin{split}
\int\limits_{\textit{D} \times S}\bar{I}^{2}(\hat{T} , s,\boldsymbol{\Omega})ds  d\boldsymbol{\Omega}   &\leq \int\limits_{\textit{D}\times S}\bar{I}^{2}(0, s,\boldsymbol{\Omega})dx  d\boldsymbol{\Omega}  + 2\nu \int\limits_{0}^{\hat{T}}\int\limits_{\textit{D}\times S }\bar{I}^{2}(t, s,\boldsymbol{\Omega})dtds  d\boldsymbol{\Omega} \\
&\quad+\nu \int\limits_{\beta}\vert \Omega.k \vert \bar{I}^{2}(t, s,\boldsymbol{\Omega})dtdk(s) d\boldsymbol{\Omega}  +U+ \int\limits_{\textit{D} \times S} \vert f\vert^{2}    d\boldsymbol{\Omega}  ds+\nu \int\limits_{\boldsymbol{D}}\frac{(\mathscr{E}(I^{\ast},\Phi))^{2}}{2}dX. \label{eq:lettaa}
\end{split}
\end{dmath}
Let
\begin{dmath}
U =2\nu \int\limits_{0}^{\bar{T}}  \int\limits_{\textit{D}\times S} \frac{k_{s}}{4\pi}\int\limits_{4\pi}\Phi(\boldsymbol{\Omega},\boldsymbol{\Omega}') \bar{I}(t, s,\boldsymbol{\Omega})\bar{I}((t, s,\boldsymbol{\Omega})d\boldsymbol{\Omega} d{\Omega'}dsdt.
\end{dmath}
The value of $U$ in Eq.(\ref{eq:lettaa}) can be determined through repeated application of the Cauchy-Schwartz inequality as follows:\\
\begin{dmath}
U \leq 2\nu \int\limits_{0}^{\bar{T}}\int\limits_{\textit{D}\times S }\bar{I}^{2}(t, s,\boldsymbol{\Omega})d{\boldsymbol{\Omega} }d{\Omega'}dsdt =\frac{2\nu(\Vert k_{s}\Vert_{L^{\infty}}+ \Vert \Sigma_{g}\Vert_{L^{\infty}})}{4\pi}\int\limits_{0}^{\bar{T}}\int\limits_{\textit{D}\times S}\bar{I}^{2}(t, s,\boldsymbol{\Omega})d{\Omega}d{\Omega'}dsdt=V_{1}\int\limits_{0}^{\bar{T}}\int\limits_{\textit{D}\times S}\bar{I}^{2}(t, s,\boldsymbol{\Omega})d{\boldsymbol{\Omega}}d{\Omega'}dsdt.~~~ 
\end{dmath}
We get after identifying the constant from
\begin{dmath}
\begin{split}
\int\limits_{\textit{D}\times S}\bar{I}^{2}(t, s,\boldsymbol{\Omega})dsd{\boldsymbol{\Omega} } &\leq \int\limits_{\textit{D}\times S}(\mathrm{R}^{\ast}_{tb})^{2}ds d\boldsymbol{\Omega} +\nu \int\limits_{\beta}(\mathrm{R}^{\ast}_{sb})^{2}dtdk(s)  d\boldsymbol{\Omega} \\
& \quad +c\Vert(\mathrm{R}^{\ast}_{int}) \Vert^{2}_{L^{2}(\textit{D}_{T} \times S)}+\nu\int\limits_{\boldsymbol{D}}(\mathscr{E}(I^{\ast},\Phi))^{2} dX+\nu V_{1}  \int\limits_{0}^{\bar{T}}\int\limits_{\textit{D}\times S }\bar{I}^{2}(t, s,\boldsymbol{\Omega})d{\boldsymbol{\Omega} }d{\Omega'}dsdt,\label{eq:lettaaa}
\end{split}
\end{dmath}
We apply the integral in Gronwall inequality form to Eq.(\ref{eq:lettaaa}) to obtain:
\begin{dmath}
\begin{split}
\int\limits_{\textit{D}\times S}\bar{I}^{2}(t,s,\boldsymbol{\Omega})dsd{\boldsymbol{\Omega} } &\leq \left( 1+\nu V_{1}\hat{T}^{2}\exp^{\nu V_{1}\hat{T}} \right)\\ 
&\quad \times \left( \int\limits_{\textit{D} \times S}(\mathrm{R}^{\ast}_{tb})^{2}ds d\boldsymbol{\Omega} +\nu \int\limits_{\beta}\mathrm({R}^{\ast}_{sb})^{2}dtdk(s)  d\boldsymbol{\Omega}  \right)\\
& \quad+  \left( 1+\nu V_{1}\bar{T}^{2}\exp^{\nu V_{1}T} \right)\\
&\quad
\times \left( v\Vert(\mathrm{R}^{\ast}_{int}) \Vert^{2}_{L^{2}(\textit{D}_{T} \times S)}+(\mathscr{E}(I^{\ast},\Phi))^{2}dX \right), \label{eq:lettaaaa}
\end{split}
\end{dmath}
Integrating Eq.(\ref{eq:lettaaaa}) over $(0,T)$
\begin{dmath}\label{eq:lettaaaaa}
\begin{split}
(\mathrm{E}_{G})^{2} = \int\limits_{\boldsymbol{D}} \bar{I}^{2}(t,s,\boldsymbol{\Omega}) \, ds \, d\boldsymbol{\Omega}  
&\leq \left( T + \nu V_{1} \hat{T}^{2} \exp^{\nu V_{1} T} \right) \left( \int\limits_{\textit{D} \times S} (\mathrm{R}^{\ast}_{tb})^2 \, ds \, d\boldsymbol{\Omega}  
+ \nu \int\limits_{\beta} (\mathrm{R}^{\ast}_{sb})^2 \, dt \, dk(s) \, d\boldsymbol{\Omega}  \right)\\
& \quad + \left( \hat{T} + \nu V_{1} \hat{T}^{2} \exp^{\nu V_{1} \hat{T}} \right) \left( v \Vert \mathrm{R}^{\ast}_{int} \Vert^{2}_{L^{2}(\textit{D}_{T} \times S)}
+ \mathscr{E}(I^{\ast}, \Phi)^{2} \, dX \right).
\end{split}
\end{dmath}
The points used for training in $\zeta$ are Sobol quadrature points, therefore the training error $E_{T}$ represents the quasi-Monte Carlo quadrature of the integral Eq.(\ref{eq:lettaaaaa}). Consequently, this aligns with the Koksma-Hlawaka inequality \cite{caflisch1998monte},\\
\begin{dmath} 
\int\limits_{\textit{D} \times S}(\mathrm{R}^{\ast}_{tb})^{2}ds d\boldsymbol{\Omega} 
\leq (\mathrm{E}_{T}^{tb})^{2}+ HK_{V}(\mathrm{R}^{\ast}_{tb})^{2}))\frac{(\log(N_{tb}))^{2d}}{N_{tb}},
\end{dmath}
Similarly
\begin{dmath}
\int\limits_{\beta} (\mathrm{R}^{\ast}_{sb})^{2} \, dt \, dk(s) \, d\boldsymbol{\Omega}  
\leq (\mathrm{E}_{T}^{sb})^{2} + HK_{V}((\mathrm{R}^{\ast}_{sb})^{2}) \frac{(\log(N_{sb}))^{2d}}{N_{sb}},
\end{dmath}
\begin{dmath}
c\Vert(\mathrm{R}^{\ast}_{int}) \Vert^{2}_{L^{2}(\textit{D}_{T} \times S)} \leq (\mathrm{E}_{T}^{int})^{2}+ HK_{V}(\mathrm{R}^{\ast}_{int})^{2})\frac{(\log(N_{int}))^{2d+1}}{N_{int}},
\end{dmath}
In this context, $(\Omega^{\boldsymbol{S}})$ for $1\leq i \leq N_{\boldsymbol{S}}$ represent points and weights of the Gauss quadrature rule $ a=a(p)$ as
\begin{dmath}
\int\limits_{\boldsymbol{D}}(\mathscr{E}(I^{\ast},\Phi))^{2}dX  \leq \bar{V}N_{\boldsymbol{S}}^{-2a}
\end{dmath}
where $\bar{V}=(\Vert\boldsymbol{D}\Vert, \Vert \Phi \Vert_{C^{p}}, \Vert \hat{I}\Vert)$.\\
\end{proof}
\begin{appendixthm}\label{thm:generalization_error_forward_inverse}
Let \( I \in L^{2}(\textit{D} \times S) \). Consider that $I$ is a unique weak solution of the RTE Eq.(\ref{eqn:BB}) with coefficients \( k_{e}, k_{s} \in L^\infty(\textit{D}) \), where \( 0 \leq k_{e}, k_{s} \), and the bounds \( \min(k_e), \max(k_e) \) and \( \min(k_s), \max(k_s) \) correspond to the lower and upper limits of these coefficients. The phase function \( \Phi(\boldsymbol{\Omega}, \boldsymbol{\Omega}') \) is assumed to lie in the space \( C^p(S \times S) \) for some \( p > 0 \).
Let \( I^* = I_{\Theta} \in C^p(\textit{D} \times S) \) represent the output of Algorithm \ref{alg1}, which addresses the forward problem for the RTE with a graded index medium. Suppose that both the Hardy-Krause variation \( V_{HK}(I^*) \) and \( V_{HK}(\mathrm{R}_{int, \Theta^*}) \) are finite, and that Sobol points \( \zeta_{int}, \zeta_{sb} \) are used as training points as previously defined. Furthermore, assume that a Gauss-quadrature rule of order \( a = a(p) \) is applied to approximate the scattering term in residual.
Now, assume \( M > 0 \), refractive index  \( n=f(s) \in D\) and that the following inequalities hold:
\[
n^{2}k_{a}I_{b}(T_{g}) \leq n^{2}k_{a}MI + f(s),
\]
and 
\[
( \boldsymbol{\Omega} \cdot \nabla)I+ \frac{1}{n \sin \theta} \frac{\partial}{\partial \phi} \left( (s_1 \cdot \nabla n) I \right) + \frac{1}{n \sin \theta} \frac{\partial}{\partial \theta} \left( I(\boldsymbol{\Omega} \cos \theta - k) \cdot \nabla n \right) \leq nMI + ( \boldsymbol{\Omega} \cdot \nabla_{s})I, \quad M > 0.
\]
Under the assumption $l>0$:
\begin{equation}\label{assumptionA}
\Big(\min(k_{e}) + \min(M) - \min(k_{s}) - 1  \Big) - \frac{2\nu(\max(k_{s}) + \| \Sigma_{g} \|_{L^{\infty}})}{4\pi}>l,
\end{equation}
In the above inequalities for \( M > 0 \), the generalization error for forward  problems is estimated by:
\begin{align}
(\mathrm{E}_{G_{steady_{f}}})^{2} &\leq V \left( \nu (\mathrm{E}_{T}^{sb})^{2} + \nu (\mathrm{E}_{T}^{int})^{2} \right) \notag \\
&+ V \left( \dfrac{(\log(N_{sb}))^{2d}}{N_{sb}} + \nu \dfrac{(\log(N_{int}))^{2d}}{N_{int}} + \nu N_{\boldsymbol{S}}^{-2a}\right),
\end{align}
Where  \( \bar{V} = (\Vert \boldsymbol{D} \Vert, \Vert \Phi \Vert_{C^p}, \Vert \hat{I} \Vert) \), \( C^{\epsilon} \) is a constant dependent on \( l \)
And
\[
V = \max \left( \frac{2}{l}, \frac{2}{l}V_{HK} \left(  R_{sb}^\ast \right)^2,  \frac{2C^{\epsilon}}{l}\left( R_{int}^\ast \right)^2,  \frac{2C^{\epsilon}}{l} \bar{V}N_{\boldsymbol{S}}^{-2a}\right).
\]
\end{appendixthm}
\begin{proof}
We are following \cite{mishra2021physics}
\begin{dmath}
\mathscr{E}_{S_{f}}(I^{\ast}, \Phi) = \sum\limits_{i=1}^{N_{\boldsymbol{S}}}w_{i}^{\boldsymbol{S}}\Phi(\boldsymbol{\Omega},\boldsymbol{\Omega}_{i}^{\boldsymbol{S}}) I^{\ast}(s,\boldsymbol{\Omega}_{i}^{\boldsymbol{S}}) -\int\limits_{4\pi}\Phi(\boldsymbol{\Omega},\boldsymbol{\Omega}') I^{\ast}(s,\boldsymbol{\Omega}) d\Omega',
\end{dmath}
$\bar{I}=I^{\ast}-I, $ 
\begin{dmath}
 (k_{e} + \boldsymbol{\Omega} \cdot \nabla)\bar{I} + \frac{1}{n\sin \theta}\frac{\partial}{\partial \theta}\left\lbrace \bar{I}(\boldsymbol{\Omega} \cos \theta - k) \cdot \nabla n \right\rbrace 
 +  \frac{1}{n\sin \theta}\frac{\partial}{\partial \phi}\left\lbrace (s_{1} \cdot \nabla n)\bar{I}\right\rbrace= n^{2}k_{a}\bar{I}_{b}(T_{g}) +\frac{k_{s}}{4\pi}\int\limits_{4\pi}\Phi(\boldsymbol{\Omega},\boldsymbol{\Omega}') \bar{I}(s,\boldsymbol{\Omega}) d\Omega'+ \mathscr{E}_{S_{f}}(I^{\ast}, \Phi), \label{eq:eq1ss} \\
\end{dmath}
\begin{align*}
\bar{I}(t, s,\boldsymbol{\Omega})&=\mathrm{R}^{\ast}_{sb}, ~~(t, s,\boldsymbol{\Omega}) \in \beta_{o},  \label{eq:eq3s}\tag{A}\\
\end{align*}
Multiplying  Eq.\eqref{eq:eq1ss} by $\bar{I}$ on both side.
\begin{dmath}
 k_{e} \bar{I^{2}} +\frac{1}{n\sin \theta}\frac{\partial}{\partial \theta}\left\lbrace \bar{I}(\boldsymbol{\Omega} \cos \theta - k) \cdot \nabla n \right\rbrace \bar{I}
 -  \frac{1}{n\sin \theta}\frac{\partial}{\partial \phi}\left\lbrace (s_{1} \cdot \nabla n)\bar{I}\right\rbrace \bar{I} =-( \boldsymbol{\Omega} \cdot \nabla_{s})\left( \frac{\bar{I}^{2}}{2}\right) + n^{2}k_{a}\bar{I}_{b}(T_{g})\bar{I}+\frac{k_{s}}{4\pi}\int\limits_{4\pi}\Phi(\boldsymbol{\Omega},\boldsymbol{\Omega}') \bar{I}( s,\boldsymbol{\Omega})\bar{I}( s,\boldsymbol{\Omega'})d{\Omega'}+ 
\mathscr{E}_{S_{f}}(I^{\ast}, \Phi)\bar{I},
\end{dmath}
We can observe that  $  \frac{1}{n\sin \theta}\frac{\partial}{\partial \phi}\left\lbrace (s_{1} \cdot \nabla n)\bar{I}\right\rbrace  + \frac{1}{n\sin \theta}\frac{\partial}{\partial \theta}\left\lbrace \bar{I}(\boldsymbol{\Omega} \cos \theta - k) \cdot \nabla n \right\rbrace  \leq M\bar{I},$\\
and $  n^{2}k_{a}\bar{I}_{b}(T_{g}) \leq n^{2}k_{a}M\bar{I} + f(s)$.\\
Integrating over $\textit{D}\times S$ and using Cauchy inequality for $k_{e}, k_{s}>0,$ 
\begin{dmath}\label{eq:callU}
\begin{split}
\int\limits_{\textit{D}\times S}\bar{I}^{2}( s,\boldsymbol{\Omega})ds  d\boldsymbol{\Omega}  &\leq  - \int\limits_{(\partial{\textit{D}} \times S)_{-}}(\Omega \cdot k(x)) \frac{\bar{I}^{2}(s,\boldsymbol{\Omega})}{2}dk(s) d\boldsymbol{\Omega}\\
& \quad + \int\limits_{\textit{D}\times S} \frac{k_{s}}{4\pi}\int\limits_{4\pi}\Phi(\boldsymbol{\Omega},\boldsymbol{\Omega}') \bar{I}(s,\boldsymbol{\Omega})\bar{I}(s,\boldsymbol{\Omega}')d\Omega' d\boldsymbol{\Omega}  ds\\ 
& \quad + \int\limits_{\textit{D}\times S} f(s) \bar{I}(s,\boldsymbol{\Omega})d\boldsymbol{\Omega}  ds \\
& \quad + \int\limits_{\textit{D}\times S}(\mathscr{E}_{S_{f}}(I^{\ast},\Phi)\bar{I}(s,\boldsymbol{\Omega}))d\boldsymbol{\Omega}  ds,
\end{split}
\end{dmath}
In this case, \( dk(s) \) signifies the surface measure on \( \partial D \), and it can defined as
\[\beta_{o}=(\partial{\textit{D}} \times S)_{-}=\left\lbrace (s,\boldsymbol{\Omega}) \in \partial{D} \times S :\boldsymbol{\Omega}.k(s)\leq 0 \right\rbrace \]  and with  $k(s)$ unit normal at  $s \in \partial{\textit{D}},$\\
Let
\begin{dmath}
U =2\nu  \int\limits_{\textit{D}\times S} \frac{k_{s}}{4\pi}\int\limits_{4\pi}\Phi(\boldsymbol{\Omega},\boldsymbol{\Omega}') \bar{I}(s,\boldsymbol{\Omega})\bar{I}((s,\boldsymbol{\Omega})d{\boldsymbol{\Omega}}d{\Omega'}ds,
\end{dmath}
Now $U$ in Eq.(\ref{eq:callU}) can be estimated by successive application of Cauchy-Schwartz inequality as
\begin{dmath}
U \leq 2\nu \int\limits_{\textit{D}\times S}\bar{I}^{2}(s,\boldsymbol{\Omega})d{\boldsymbol{\Omega}}d{\Omega'}ds =\frac{2\nu(\Vert k_{s}\Vert_{L^{\infty}}+ \Vert \Sigma_{g}\Vert_{L^{\infty}})}{4\pi}\int\limits_{\textit{D}\times S}\bar{I}^{2}( s,\boldsymbol{\Omega})d{\boldsymbol{\Omega}}d{\Omega'}ds=V_{1}\int\limits_{\textit{D}\times S}\bar{I}^{2}(s,\boldsymbol{\Omega})d{\boldsymbol{\Omega}}d{\Omega'}ds,~~~ 
\end{dmath}
Utilizing integration by parts, and applying the presumed limits on \( k_e \), \( k_s \), and \( M \), we derive:\\
\begin{dmath}
\begin{split}
\Big( \min(k_{e}) + \min(M) - \min(k_{s}) - 1 \Big)\int\limits_{\textit{D}\times S}\bar{I}^{2}( s,\boldsymbol{\Omega})dsd{\boldsymbol{\Omega}} & \leq \nu \int\limits_{\beta_{o}}(\mathrm{R}^{\ast}_{sb})^{2}dk(s)  d\boldsymbol{\Omega} \\
&\quad +\int\limits_{\textit{D}\times S} (\mathrm{R}^{\ast}_{int})\bar{I}d\boldsymbol{\Omega}  ds +\nu\int\limits_{\textit{D} \times S}(\mathscr{E}_{S_{f}}(I^{\ast},\Phi))\bar{I} dX\\
&+\nu V_{1}  \int\limits_{\textit{D}\times S}\bar{I}^{2}(s,\boldsymbol{\Omega})d{\boldsymbol{\Omega}}d{\Omega'}ds,
\end{split}
\end{dmath}
Given assumption \ref{assumptionA}, $\exists$ an $\epsilon > 0$ such that

\[\Big( \min(k_{e}) + \min(M) - \min(k_{s}) - 1 \Big)- \frac{2\nu(\max(k_{s})+ \Vert \Sigma_{g}\Vert_{L^{\infty}})}{4\pi}  - 2\epsilon > \frac{l}{2},
\]
We apply the $\epsilon$-version of Cauchy's inequality.
\begin{dmath}
\int\limits_{\textit{D} \times S} \bar{I}^{2}(s, \boldsymbol{\Omega}) \, ds \, d\boldsymbol{\Omega} 
\leq 
\frac{2}{l}\int\limits_{\beta_{o}} \left( \mathrm{R}^{\ast}_{sb} \right)^{2} \, dk(s) \, d\boldsymbol{\Omega} 
+\frac{2C^{\epsilon}}{l}  \int\limits_{\textit{D} \times S} \left( \mathrm{R}^{\ast}_{int} \right)^{2} \, d\boldsymbol{\Omega} 
+ \mathscr{E}_{S_{f}}\left( I^{\ast}, \Phi \right)^{2} \, dX.
\end{dmath}
\end{proof}
\bibliographystyle{abbrv}
\bibliography{sample}
\end{document}